\documentclass[11pt]{article}
\usepackage{color}

\usepackage{amsfonts}
\usepackage{amscd}
\usepackage{amstext}
\usepackage{latexsym}
\usepackage{amsmath}
\usepackage{empheq}
\usepackage{amsbsy}
\usepackage{graphics}
\usepackage{color}
\usepackage[english]{babel}
\usepackage[latin1]{inputenc}
\numberwithin{equation}{section}
\usepackage{times}

\newtheorem{theorem}{Theorem}[section]
\newtheorem{lemma}{Lemma}[section]
\newtheorem{proposition}{Proposition}[section]
\newtheorem{corollary}{Corollary}[section]

\newcommand{\cqd}{\hfill $\Box$ \vspace{0.5cm}}
\newcommand{\rn} {\mathbb{R}^N}
\newcommand{\dd}{\,\mathrm{d}}
\newcommand{\cinfsc}{C_c^\infty}
\newcommand{\grad}{\nabla}

\newcommand{\ol}{{\Omega_\lambda}}

\newcommand{\al}{{A_\lambda}}

\newcommand{\dem}{\flushleft{{\emph{Proof.}}}}
\newcommand{\imerso}{\hookrightarrow}
\newcommand{\tofraco}{\rightharpoonup}
\newcommand{\nat}{n \in \mathbb N}

\newcommand{\re}{\mathbb R}

\newcommand{\ci}{c_{\infty}}

\newcommand{\supp}{\textrm{supp}}

\title{Nontrivial solutions for a mixed boundary problem for  Schr\"{o}dinger equations with an external magnetic field\thanks{C. O. Alves and S. H. M. Soares were partially supported by CNPq. R. C. M. Nemer was supported by   FAPESP.}}

\author{Claudianor O. Alves \\
\noindent Unidade Acad\^emica de Matem\'atica e Estat\'istica\\
\noindent Universidade Federal de Campina Grande \\
\noindent 58109-970, Campina Grande - PB, Brazil.\\
\noindent e-mail: {\tt{coalves@dme.ufcg.edu.br}}\\ \\
\noindent Rodrigo C. M. Nemer\ \ and\ \ S\'{e}rgio H. M. Soares\thanks{Corresponding author: phone 55 (16) 3373-9660; fax 55 (16) 3373-9650.} \\
\noindent Departamento de Matem\'atica \\
\noindent Instituto de Ci\^{e}ncias Matem\'{a}ticas e de
 Computa\c{c}\~{a}o\\
\noindent Universidade de S\~ao Paulo \\
\noindent 13560-970, S\~ao Carlos - SP, Brazil. \\
\noindent e-mail: {\tt{rcmnemer@icmc.usp.br, monari@icmc.usp.br}}\\
}
\date{}

\begin{document}

\maketitle

\begin{abstract}
We study the existence of solutions for a class of nonlinear Schr\"{o}dinger equations involving a magnetic field with mixed
Dirichlet--Neumann boundary conditions.  We use Lyusternik-Shnirelman category and the Morse theory to estimate the number of nontrivial solutions in terms of the topology of the part of the boundary where the Neumann condition is prescribed.
\end{abstract}



\section{Introduction}

A major role in quantum physics is played by the nonlinear Schr\"odinger equation
\begin{equation}\label{1}
ih\frac{\partial \Psi}{\partial t}= \left( \frac{h}{i}\nabla - A(x)\right)^2\Psi + U(x)\Psi - f(|\Psi|^2)\Psi, \quad x \in \Omega,
\end{equation}
where $\Omega$ is a bounded smooth domain in $\mathbb{R}^N$, $N\geq 3$,  $t \in \mathbb{R}$, $h$ is a positive constant, $i$ is the imaginary unit, $\Psi :  \mathbb{R} \times \mathbb{R}^N \to \mathbb{C}$ is the wave function, $f$ is a  nonlinear term, $U$ is the real electric potential,  $A : \mathbb{R}^N \to \mathbb{R}^N$ denotes a magnetic potential and the Schr\"odinger operator is defined by
\[
\left( \frac{h}{i}\nabla - A(x)\right)^2\Psi = -h^2\Delta \Psi - \frac{2h}{i}A\nabla \Psi + |A|^2\Psi - \frac{h}{i}\Psi {\rm{div}} A.
\]
 We are interested in standing wave solutions, that is, solutions for (\ref{1}) in the form $\Psi (t,x) = e^{-iEt/h}u(x)$, where $u$ satisfies
\begin{equation}\label{2}
\left( \frac{h}{i}\nabla - A(x)\right)^2u + V(x)u =  f(|u|^2)u, \quad x \in \Omega,
\end{equation}
where $V (x) = U(x) - E$.  Assuming that $V\equiv 1$, it follows immediately that $u$ is a solution of \eqref{2} if, and only if, the function $v(x) = u(hx)$ solves
\begin{equation}\label{2'}
\left( \frac{1}{i}\nabla - A_\lambda(x)\right)^2v +  v =  f(|v|^2)v, \quad x \in \Omega_\lambda,
\end{equation}
where $\lambda = h^{-1}$, $A_\lambda(x) = A(\lambda^{-1}x)$ and $\Omega_\lambda \doteq \lambda \Omega$. The case with no magnetic vector field, namely $A=0$, has been widely studied in the literature. We refer to  \cite{ABC},  \cite{ambrosetti-ruiz}, \cite{benci-cerami}, \cite{delpino-felmer},  \cite{floer-wiestein}, \cite{rabinowitz}, \cite{wang},  \cite{wang-zeng},  and references in these papers. Existence results for the magnetic case were established in \cite{alves-figueiredo-furtado}, \cite{chabrowski-szulkin},  \cite{cingolani-clapp}, \cite{cingolani-jeanjean-secchi}, \cite{cingolani-secchi2005}, \cite{cingolaniJDE2003}, \cite{cingolani-secchi2002},  \cite{esteban-lions},  \cite{kurata}, \cite{tang}.  In \cite{alves-figueiredo-furtado}, the authors have proved that
if  $f$ is a superlinear function with subcritical growth, then for large values of $\lambda>0$, the equation \eqref{2'} with boundary Dirichlet condition has at least ${\rm{cat}}_{\Omega_\lambda}(\Omega_\lambda)$ nontrivial weak solutions, where  ${\rm{cat}}_{\Omega_\lambda}(\Omega_\lambda)$ denotes the
 the Ljusternik-Schnirelman category of $\Omega_\lambda$ in $\Omega_\lambda$. In the seminal work \cite{benci-cerami}, Benci and  Cerami used Ljusternik-Schnirelman category and  Morse theory to estimate the number of positive solutions of the problem
\begin{equation}
  \left\{
\begin{array}{ll}
-\epsilon \Delta u + u = f(u), & \mbox{in $\Omega$,}\\
\quad u= 0, &  \mbox{on $\Omega$,}
\end{array}
\right. \label{benci-cerami}
\end{equation}
where $\Omega$ is a bounded domain. It is proved  that for $\epsilon$ sufficiently small the number of positive solutions is at least ${\rm{cat}}_{\Omega}(\Omega)$. They also proved via Morse theory that the number of solutions depends on the topology of $\Omega$, actually on  $\mathcal{P}_t(\Omega)$,  the Poincar\'{e} polynomial of $\Omega$.  In \cite{candela-lazzo}, Candela and Lazzo have considered this same equation with mixed Dirichlet-Neumann boundary conditions with $f(t)=|t|^{p-2}t$. It was proved that the number of positive solutions is influenced by the topology of the part $\Gamma_1$ of the boundary $\partial\Omega$ where the Neumann condition is assumed, more precisely, if $(N-1)$-dimensional Lebesgue measure in $\mathbb{R}^N$ is positive, then the respective problem has at least category of a set $\Gamma_1$, provided $\epsilon$ is sufficiently small.

Motivated by the results just described, a natural question is whether same kind of result holds for the mixed boundary problem with magnetic field
  \begin{equation}
\left\{\begin{array}{rcll}
            \left( -i\nabla - A_\lambda \right)^2 u + u & = & f(|u|^2) u,& \textrm{ in } \Omega_\lambda \\
            u & = & 0,  &\textrm{ on } \Gamma_{0 \lambda}\\
            \frac{\partial u}{\partial \nu} & = & 0, & \textrm{ on } \Gamma_{1 \lambda}
        \end{array}
\right. \label{2pal}
\end{equation}
where $\lambda$  is a positive real parameter,  $\Omega_\lambda = \lambda \Omega$ is an expanding set,  $\Omega \subset \mathbb R ^N$ ($N \ge 3$) is a bounded domain with smooth boundary $\partial \Omega = \overline{\Gamma_0} \cup \overline{\Gamma_1}$, where $\Gamma_0, \Gamma_1$ are smooth disjoint submanifolds with positive $(N-1)$-dimensional Lebesgue measure in $\mathbb{R}^N$,  $\Gamma_{0 \lambda} \doteq \lambda \Gamma_0$,  $\Gamma_{1 \lambda} \doteq \lambda \Gamma_1$, $A \in C( \Omega,\mathbb{R}^N)$ and $f\in C^1(\mathbb{R^+})$ satisfies:
\begin{itemize}
\item[($f_{1})$] $f(s) = o(1)$  and $f'(s) = o(1/s)$, as $s\to 0^+$.
\item[($f_{2})$] There exists  $q \in (2, 2^*)$ such that
\[
\lim_{s \to \infty}\displaystyle\frac{f(s)}{s^{\frac{q-2}{2}}} = 0\quad \mbox{and}\quad \lim_{s \to \infty}\displaystyle\frac{f'(s)}{s^{\frac{q-4}{2}}} = 0,
\]
where $2^* = {2N}/{(N-2)}$.
\item[($f_{3})$] There exists $\theta > 2$ such that
\[
0 < \frac{\theta}{2} F(s) \le sf(s), \textrm{ for } s>0,
\]
where $F(s) = \displaystyle\int_0^s f(t) \dd t$.
\item[($f_{4})$] $f'(s) > 0$, for all $s > 0$.
\item[($f_{5}$)] There exist $q \in (2, 2^*)$ and a constant $C>0$ such that
\[
sf(s) - F(s) \ge C|s|^{q/2}, \textrm{ for all } s \ge 0.
\]
\end{itemize}

We state that the  magnetic field does not play any role on the number of solutions of (\ref{2}) and therefore a result in the same spirit of  \cite{benci-cerami} and \cite{candela-lazzo} holds. More precisely, our main results are the following:

\begin{theorem} \label{2t1}
Suppose that $f$ satisfies $(f_1)-(f_5)$. There exist $\lambda^* > 0$ such that for any $\lambda > \lambda^*$  problem \eqref{2pal} has at least $cat_{\Gamma_{1\lambda}}(\Gamma_{1\lambda})$  nontrivial  weak solutions.
\end{theorem}

To established the result in terms of Morse theory, we introduce some notation. For any $\lambda > 0$,  let  $H^1_{A_\lambda}(\ol, \Gamma_{0 \lambda})$ be the Hilbert space
\[
H^1_{A_\lambda}(\ol, \Gamma_{0 \lambda}) \doteq \{ u \in L^2(\ol, \mathbb C); |\grad_{\al} u| \in L^2(\ol), \textrm{trace of } u = 0 \textrm{ on } \Gamma_{0 \lambda} \},
\]
endowed with the norm
\[
<u, v>_\al \doteq \textrm{Re} \left\{ \int_{\ol} (\grad_\al u \overline{\grad_\al v} + u \overline v)  \dd x \right\},
\]
where
\[
\grad_{A_\lambda} u \doteq (D^j_{A_\lambda} u)^N_{j = 1}, \quad D^j_{A_\lambda} u \doteq -i \partial_j u - A^j_\lambda u
\]
and $\textrm{Re}(w)$  is the real part of $w \in \mathbb C$ and  $\overline w$  is its complex conjugate. The norm induced by this inner product is given by
\[
\| u \|_{A_\lambda} = \left( \int_{\ol}(|\grad_\al u |^2 + |u|^2)  \dd x  \right)^{1/2}.
\]
By \cite{esteban-lions},  we can state a version of diamagnetic
inequality for the space $H^1_{\al} (\ol, \Gamma_{0 \lambda})$: For any $u\in H^1_{\al} (\ol, \Gamma_{0 \lambda})$,
\begin{equation}\label{diamagnetic}
|\grad_{\al}u| \ge |\grad|u||.
\end{equation}
As a consequence, the embedding $H^1_{\al} (\ol, \Gamma_{0 \lambda}) \imerso L^p(\ol, \re)$ is continuous for  $1 \leq p \leq 2^*$ and it is compact for $1\leq p < 2^{*}$. It is worth pointing out that the embedding constants do not depend on $\lambda$, because of the assumption that $\Omega \subset \mathbb R ^N$ ($N \ge 3$) is a bounded domain with smooth boundary $\partial \Omega$. {We also emphasize that the regularity on $\partial \Omega$ assumed here must be sufficient to obtain  $r_0 > 0$  such that
\begin{equation}\label{CEEIU}
 B_{r_0}(y + r_0 \nu_y) \subset \Omega\quad \mbox{and}\quad
 B_{r_0}(y - r_0 \nu_y) \subset \rn \backslash \Omega,
\end{equation}
uniformly for  $y \in \partial\Omega$, where $\nu_y$ is the inward unitary normal vector to $\partial \Omega$ in $y$ and $B_r(z)$ denotes the ball of radius $r$ centered at $z$.}

The functional associated with \eqref{2pal} $I_{\lambda}: H^1_{A_\lambda}(\ol, \Gamma_{0 \lambda}) \to \re$ is given by
\begin{equation}\label{functional}
I_{\lambda} (u) \doteq \frac{1}{2} \int_{\ol}|\grad_{A_\lambda} u|^2 + |u|^2 \dd x - \frac{1}{2}\int_{\ol} F(|u|^2) \dd x.
\end{equation}
From $(f_1)-(f_2)$, the functional $I_{\lambda}$ is well defined and belongs to  $C^2(H^1_{A_\lambda}(\ol, \Gamma_{0 \lambda}), \mathbb{R})$. Furthermore,
\[
I'_{\lambda} (u)v \doteq \textrm{Re} \left\{\int_{\ol} \grad_{A_\lambda} u \overline{\grad_{A_\lambda} v} + u \overline{v} \dd x - \int_{\ol} f(|u|^2)u\overline{v} \dd x \right\},
\]
for all $u, v \in H^1_{A_\lambda}(\ol, \Gamma_{0 \lambda})$. Thus, every critical point of $I_{\lambda}$ is a weak solution of \eqref{2pal}.

In the notation of \cite{benci-cerami}, we have if $u$ is an isolated critical point of $I_{\lambda}$ and $I_{\lambda}(u)=c$, the polynomial Morse index  $i_t(u)$ of $u$ is defined by
\[
i_t(u) = \sum_k \text{dim}[H^k(I_{\lambda}^c\cap U, (I_{\lambda}^c\setminus \{u\})\cap U)]t^k,
\]
where $ H^k(\cdot, \cdot)$ denotes the $k$th group de homology with coefficients in some field $\mathbb{K}$, $U$ is a neighborhood of $u$ and $$I_{\lambda}^c = \{ v \in H^1_{A_\lambda}(\ol, \Gamma_{0 \lambda}); I_{\lambda}(v) \leq c\}.
$$

As is proved in \cite[Theorem I.5.8]{benci}, if $u$ is a non-degenerate critical point, then $i_t(u) = t^{\mu(t)}$, where $\mu(u)$ denotes the numeric Morse index of $u$.

Let $X$ be a topological space. The Poincar\'{e} polynomial of $X$ is defined by
\[
\mathcal P_t(X) = \sum_k \text{dim}[H^k(X)]t^k.
\]

Following \cite{benci-cerami}, we can prove the ensuing multiplicity result:
\begin{theorem} \label{2t2}
Suppose that $f$ satisfies $(f_1)-(f_5)$ and the set $\mathcal K$  of nontrivial solutions of problem \eqref{2pal} is discrete. Then, there exists
 $\lambda^* > 0$ such that
 \[
\sum_{u \in \mathcal K} i_t(u) = t \mathcal P_t(\Gamma_{1\lambda}) + t^2 [\mathcal P_t (\Gamma_{1\lambda}) - 1] + (t+1)\mathcal Q (t),
\]
 for every $\lambda > \lambda^*$, where $\mathcal Q(t)$  is a polynomial with non-negative integer coefficients.
\end{theorem}

In the non-degenerate case, we have:
\begin{corollary} \label{2c1}
Suppose that $f$ satisfies $(f_1)-(f_5)$ and the solutions of problem \eqref{2pal} are non-degenerate. Then, there exists
 $\lambda^* > 0$ such that
\[
\sum_{u \in \mathcal K} t^{\mu(u)} = t \mathcal P_t(\Gamma_{1\lambda}) + t^2 [\mathcal P_t (\Gamma_{1\lambda}) - 1] + (t+1) \mathcal Q (t),
\]
for every $\lambda > \lambda^*$, where $\mathcal Q(t)$ is a polynomial with non-negative integer coefficients.
\end{corollary}

As observed in \cite{benci-cerami} (see also \cite{cingolani-vannella}), the application of the Morse theory can give better information than the use of the Ljusternik-Schnirelman theorem. Theorem \ref{2t2} shows that the problem  \eqref{2pal} possesses at least  $2\mathcal{P}_1(\Gamma_{1\lambda}) -1$ {nontrivial weak solutions}. In the case of $\Gamma_{1\lambda}$ is topologically trivial, we have $\mathcal{P}_1(\Gamma_{1\lambda}) =1$ and this theorem does not provide  any additional information about multiplicity of solutions. On the other hand, when $\Gamma_{1\lambda}$ is  a topologically rich  domain, for example, if $\Gamma_{1\lambda}$ is obtained by contractible submanifold cutting off $k$ contractible open non-empty sets in $\partial \Omega$, we obtain that the number of {nontrivial} solutions of \eqref{2pal} is affected by $k$, even if the category of $\Gamma_{1\lambda}$ is 2.

{In order to prove Theorems \ref{2t1} and \ref{2t2}, we combine the {Benci and Cerami}  approach \cite{benci-cerami} with a variation of the arguments of Candela and Lazzo \cite{candela-lazzo}.  The major steps in {Benci and Cerami} approach are the analysis of the behavior of some critical  levels related to problem \eqref{benci-cerami} and the comparison of the topology of $\Omega$ with  some sublevel sets of the functional associated  with \eqref{benci-cerami}. Although we use this machinery,  we have to make a detailed analysis of the behavior of the minimax levels associated with the problem \eqref{2pal} and a more involved proof that the barycenter function maps suitable sublevel sets of the functional associated with \eqref{2pal} in a  neighborhood of the portion of the boundary where the Neumann condition is prescribed. This is because  the equation \eqref{2pal} involves a magnetic field and mixed Dirichlet-Neumann boundary conditions. Moreover,  as the nonlinearity is not  necessarily homogeneous, our arguments are different from what can already be found in \cite{candela-lazzo}.  Once these crucial steps are verified,  we can employ the Morse theory developed in \cite[Section 5]{benci-cerami}  to estimate the number of nontrivial solutions to \eqref{2pal} in terms of the topology of the part of the boundary where the Neumann condition is assumed. }

{Theorems \ref{2t1} and \ref{2t2} can be seen as a complement  of  the studies made in \cite{alves-figueiredo-furtado}, \cite{benci-cerami} and \cite{candela-lazzo} in the following aspects:  1) In \cite{alves-figueiredo-furtado} only the Dirichlet  boundary condition was considered; 2) In \cite{benci-cerami}, the problem was considered for the Laplacian operator and Dirichlet  boundary condition. Here we are working with a more general boundary condition and with a class of  operators  which includes the Laplacian operator as a particular case; 3) In \cite{candela-lazzo},  the problem was  also considered for Laplacian operator and with a homogeneous nonlinearity. In the present paper we deal with a class of  nonlinearities that has the homogeneous functions as a particular case. As we are  mainly considering a non homogeneous nonlinearity, our estimates are more delicate and we need to make a  careful analysis in several estimates involving different arguments from those used in \cite{candela-lazzo}, see Sections 3, 4 and 5. }

\section{The Palais-Smale condition}

In this section we establish the Palais-Smale condition for the functional $I_\lambda$, defined by \eqref{functional}, and for the  functional $I_{\lambda}$ constrained to ${M_{\lambda}}$. As a direct consequence of  $(f_1)-(f_3)$, we obtain
\begin{itemize}
\item[$(f_{6})$] Given $\epsilon > 0$, there exist constant $C_\epsilon > 0$ such that
\[
f(s) \le \epsilon + C_\epsilon s^{{(q-2)}/{2}},\quad \forall\, s \geq 0,
\]
where $q \in (2, 2^*)$.
\item[$(f_{7})$] There exists $\theta >2$ and a constant $C > 0$ such that
\[
F(s) \ge C|s|^{\theta/{2}} - C,\  \forall\, s \ge 0,
\]
where $F(s) = \displaystyle\int_0^s f(t) \dd t$ .
\end{itemize}

\begin{proposition}\label{ps}
The functional $I_{\lambda}$ satisfies the Palais-Smale condition, that is, every sequence
$(u_n) \subset  H^1_{A_\lambda}$ for which  $\sup_{n\in \mathbb{N}}|I_{\lambda} (u_n)|<  \infty$  and  $I'_{\lambda} (u_n) \to 0$, as $n\to \infty$,  possesses a converging subsequence.
\end{proposition}
\dem \  Given a sequence $(u_n) \subset  H^1_{A_\lambda}(\ol, \Gamma_{0 \lambda})$ such that  $\sup_{n\in \mathbb{N}}|I_{\lambda} (u_n)| < \infty$ and  $I'_{\lambda} (u_n) \to 0$, as $n\to \infty$, we may assume that $I_{\lambda} (u_n) \to d$ and  $I'_{\lambda} (u_n) \to 0$, as $n\to \infty$, for some  $d \in \mathbb{R}$.  We claim that $(u_n)$ is bounded. In fact, from   $(f_{3})$, we have
\begin{align*}
d + o_n(1) + o_n(1)\|u_n\|_{\al} & = I_{\lambda} (u_n) - \frac{1}{\theta} I_\lambda'(u_n)u_n \\
                                 & = \left(\frac{1}{2} - \frac{1}{\theta}\right) \|u_n\|^2_{\al} + \int_{\ol} \left(\frac{1}{\theta} f(|u_n|^2)|u_n|^2 - \frac{1}{2} F(|u_n|^2\right) \\
                                 & \ge \left(\frac{1}{2} - \frac{1}{\theta}\right) \|u_n\|^2_{\al},
\end{align*}
where $o_n(1)$ denotes a quantity going to zero zero as $n\to \infty$. From this, we obtain that $(u_n)$ is bounded. As a consequence, we may assume that  $(u_n)$ has a subsequence, still denoted by $(u_n)$, and there exists
 $u \in H^1_{\al}(\ol, \Gamma_{0 \lambda})$  such that
\begin{empheq}[left=\empheqlbrace]{align}
u_n \tofraco u\quad  &\textrm{in } H^1_{\al}(\ol, \Gamma_{0 \lambda}), \nonumber\\
u_n \to u \quad& \textrm{ in } L^p(\ol, \mathbb C), \forall p \in [\left. 1, 2^*)\right.,  \label{propriedade}\\
 u_n \to u\quad  & \textrm{ a.e. in } \ol. \nonumber
\end{empheq}
Invoking the definition of $I'_{\lambda}$, we obtain
\[
\|u_n - u\|^2_{\al} = (I'_{\lambda} (u_n) - I'_{\lambda} (u))(u_n - u) - \textrm{Re}\left\{ \int_{\ol} (f(|u_n|^2)u_n - f(|u|^2)u) \overline{(u_n - u)}\right\}.
\]
Thus, from $(f_6)$ and (\ref{propriedade}),
\begin{eqnarray*}
\|u_n - u\|^2_{\al}  &\le&  |I'_{\lambda} (u_n) (u_n - u)| + |I'_{\lambda} (u) (u_n - u)| \\
& &+ \int_{\ol} |f(|u_n|^2)u_n - f(|u|^2)u| |u_n - u| =o_n(1),
\end{eqnarray*}
as $n\to \infty$. Hence,  $u_n \to u$ in $H^1_{\al}(\ol, \Gamma_{0 \lambda})$. \cqd

By $(f_6)-(f_7)$, it is a simple matter to check that $I_\lambda$ satisfies the geometric hypotheses of the mountain pass theorem. From this and  Proposition \ref{ps}, for any  $\lambda >0$, there exists $u_\lambda \in H^1_{\al}(\ol, \Gamma_{0 \lambda})$ such that $I'_\lambda(u_\lambda) = 0$ and $I_\lambda(u_\lambda) \doteq b_\lambda$, where $b_\lambda$ denotes the mountain pass level of the functional $I_\lambda$. From $(f_4)$, the level $b_\lambda$ satisfies (see \cite{Willem})
\begin{equation}\label{blambda}
b_\lambda = \inf_{u \in M_{\lambda}} I_\lambda (u),
\end{equation}
where $M_{\lambda}$ denotes the Nehari manifold associated with $I_\lambda$, namely
\[
M_{\lambda} = \{ u \in H^1_{A_\lambda}(\ol, \Gamma_{0 \lambda})\backslash \{0\};\  I'_{\ol} (u)u = 0\}.
\]
Since we are intend to consider the functional $I_\lambda$ constrained to $M_{\lambda}$, the next two results are required.
\begin{proposition} \label{limnehari} Suppose that $f$ satisfies $(f_1)$ and $(f_2)$. Then, there exists  $\delta_0  > 0$ independent of $\lambda > 0$ such that every  $u \in M_{\lambda}$ satisfies
\begin{equation}
\|u\|_{\al} \ge \delta_0 \quad \textrm{and} \quad I_{\lambda}(u) \ge \delta_0.
\end{equation}
\end{proposition}
\dem  \  From $(f_{6})$, given  $\epsilon > 0$ there exists $C_{\epsilon} >0$ such that for every   $u \in M_{\lambda}$,
\[
\|u\|^2_{\al}  = \int_{\ol} f(|u|^2)|u|^2 \le \epsilon \int_{\ol}|u|^2_{\ol, 2} + C_{\epsilon}\int_{\ol} |u|^q_{\ol, q}.
\]
Since  the embedding $H^1_{\al}(\ol, \Gamma_{0 \lambda}) \imerso L^p(\ol, \mathbb C)$ is continuous for $p\in [1, 2^*]$ and the embedding constant does not depend on $\lambda$, there exists a positive constant $C$  independent of $\lambda$ such that
\begin{equation}
\|u\|^2_{\al}   \leq C (\epsilon \|u\|^2_{\al} + C_\epsilon \|u\|^q_{\al}). \label{limnehari1}
\end{equation}
Taking  $\epsilon  ={1}/(2 C)$ in \eqref{limnehari1},  we have
\begin{equation}
\label{limnehari2}
\|u\|_{\al} \ge \left(\frac{1}{2 C C{_\epsilon}}\right)^{\frac{1}{q-2}} =: \delta_1 > 0.
\end{equation}
For any $u \in M_\lambda$, from $(f_{3})$ and  \eqref{limnehari2}, it follows that
\begin{eqnarray*}
I_{\lambda} (u) & = & \left(\frac{1}{2} - \frac{1}{\theta}\right)\|u\|^2_{\al} + \int_{\ol} \left( \frac{1}{\theta}f(|u|^2)|u|^2 - \frac{1}{2}F(|u|^2) \right)\\
& \ge& \left(\frac{1}{2} - \frac{1}{\theta}\right)\|u\|^2_{\al}
            \ge \left(\frac{1}{2} - \frac{1}{\theta}\right) \delta_1^2 =: \delta_2.
\end{eqnarray*}
Taking   $\delta_0 \doteq \min\{\delta_1, \delta_2\} = \delta_2$,  we conclude the proof.  \cqd
\begin{proposition} \label{psr} The functional  $I_{\lambda}$ constrained to ${M_{\lambda}}$ satisfies the Palais-Smale condition.
 \end{proposition}

\dem\  Let  $(u_n) \subset M_{\lambda}$ be a sequence such
{$\sup_{n\in \mathbb{N}} |I_\lambda(u_n)| < \infty$} and $(I_{\lambda}\big|_{M_{\lambda}})'(u_n) \to 0$, as $n\to \infty$.  We can assume, by taking a subsequence if necessary, that $ I_\lambda(u_n) \to d$, for some $d\in \mathbb{R}$.  By  \cite[Proposition 5.12]{Willem}, for each $n\in \mathbb{N}$ there exists $\mu_n \in \re$ such that
\begin{equation}
\label{psr0}
I'_{\lambda} (u_n) - \mu_n G'_{\lambda}(u_n) = (I_{\lambda}\big|_{M_{\lambda}})'(u_n) = o_n(1),
\end{equation}
where
\[
G_{\lambda}(v) = I'_{\lambda}(v)v,\quad  \forall v \in H^1_{\al} (\ol, \Gamma_{0 \lambda}).
\]
  As in the proof of  Proposition \ref{ps}, $(u_n)$  is bounded. Hence, we may suppose that $(u_n)$ has a  subsequence, still denoted by $(u_n)$, and there exists  $u \in H^1_{\al}(\ol, \Gamma_{0 \lambda})$ such that
\begin{empheq}[left=\empheqlbrace]{align}
& u_n \tofraco u \textrm{ in } H^1_{\al}(\ol, \Gamma_{0 \lambda}),  \nonumber\\
& u_n \to u \textrm{ in } L^p(\ol, \mathbb C), \forall p \in [\left. 1, 2*)\right.,  \nonumber\\
& u_n \to u \textrm{ a.e. on } \ol.  \nonumber
\end{empheq}
Since  $u_n \in M_{\lambda}$, the condition $(f_4)$ implies
\begin{equation}
G'_{\lambda} (u_n)u_n =  - 2 \int_{\ol} f'(|u_n|^2)|u_n|^4   \le  0. \label{psr1}
\end{equation}
Moreover,  by Proposition \ref{limnehari}, we have
\begin{equation}
\label{psr2}
\int_{\ol} f(|u_n|^2)|u_n|^2 \ge {\delta_0^{2}}, \quad \forall \nat.
\end{equation}
Taking $n\to \infty$ and using the Sobolev embedding, we obtain
\begin{equation}
\label{psr2.2}
\int_{\ol} f(|u|^2)|u|^2 \ge {\delta_0^{2}},
\end{equation}
and so, $ u \not\equiv 0$. From this, \eqref{psr1} and  Fatou lemma, we have
\begin{align*}
\liminf_{n \to \infty} G'_{\lambda}(u_n)u_n & = \liminf_{n \to \infty} -2 \int_{\ol} f'(|u_n|^2)|u_n|^4  \le -2 \int_{\ol } f'(| u|^2)| u|^4 < 0.
\end{align*}
Now, we use  \eqref{psr0} to obtain that  $\mu_n \to 0$, as $n\to \infty$. Consequently,  the sequence $(u_n)$ also satisfies
  $\sup_{n\in \mathbb{N}}|I_{\lambda} (u_n)|<  \infty$  and  $I'_{\lambda} (u_n) \to 0$, as $n\to \infty$.  Proposition \ref{psr} now shows that the functional  $I_{\lambda}$ constrained to ${M_{\lambda}}$ satisfies the Palais-Smale condition. \cqd

We can proceed analogously to the proof of Proposition \ref{psr} to show the next result.
\begin{corollary} \label{pcr}
I f $u $ is a critical of the functional  $I_\lambda$ constrained to  ${M_\lambda}$, then $u$ is a nontrivial critical point of  $I_\lambda$.
\end{corollary}

\section{Preliminaries}
Firstly we introduce some notation. Let $\rn_+ = \{(x_1, \dots, x_N) \in \mathbb{R}^N : x_N >0\}$ and $\mathbb R ^{N-1}= \{(x_1, \dots, x_N) \in \mathbb{R}^N : x_N = 0\}$.  Consider the problems
\begin{equation}\label{2pinf}
\left\{\begin{array}{rcl}
            -\Delta u + u & = & f(u^2) u\  \textrm{ in } \rn_+, \\
            \displaystyle\frac{\partial u}{\partial \nu} & = & 0 \textrm{ on } \mathbb R ^{N-1}
        \end{array}
\right.
\end{equation}
and
\begin{equation}\label{2prn}
\left\{\begin{array}{l}
            -\Delta u + u\  = \ f(u^2)u\quad   \textrm{ in } \rn, \\
           \  u \in H^1(\rn).
        \end{array}
\right.
\end{equation}
Consider now the respective functionals associated with the above problems
\[
J_{\infty} (u) \doteq \frac{1}{2} \int_{\rn_+} (|\grad u|^2 + u^2) - \frac{1}{2}\int_{\rn_+} F(u^2), \quad \forall u \in H^1(\rn_+),
\]
and
\[
J_{\rn} (u) \doteq \frac{1}{2} \int_{\rn} (|\grad u|^2 + u^2) - \frac{1}{2}\int_{\rn} F(u^2), \quad \forall u \in H^1(\rn).
\]
We define the corresponding Nehari manifolds and mountain pass levels:
\[
N_\infty \doteq \{ u \in H^1(\rn_+) \backslash \{0\};\  J'_\infty(u)u = 0\}
\quad \mbox{and}\quad
c_\infty \doteq \inf\limits_{N_\infty} J_{\infty},
\]
and
\[
N_{\rn} \doteq \{ u \in H^1(\rn) \backslash \{0\};\  J'_{\rn}(u)u = 0\}
\quad \mbox{and}\quad
c_{\rn} \doteq \inf\limits_{N_{\rn}} J_{\rn}.
\]
{By \cite{berestycki-lions, strauss},}  \eqref{2prn} has a radially symmetric positive solution $w \in H^1(\rn) \cap C^2(\rn)$.   Moreover, the restriction of $w$ to ${\rn_+}$ is a solution of \eqref{2pinf}. As a consequence,
\begin{equation}
\label{crncinf}
c_{\rn} = 2 c_\infty.
\end{equation}
Let  $r>0$ be such that the sets
\[
\Gamma_1^+ \doteq \{ x \in \rn ;\  \textrm{dist}(x, \Gamma_1) < r \}, \quad \Gamma_1^- \doteq \{ x \in \Gamma_1 ;\  \textrm{dist}(x, \Gamma_0) \ge r\}
\]
are homotopically equivalent to  $\Gamma_1$.  Let  $\eta \in C^{\infty}(\re_+)$ be a non-increasing function  such that
$\eta = 1$ on  $[0,r/2]$, $\eta = 0$ on $[r,+\infty)$, $|\eta'| \in L^{\infty}(\re_+)$. We will denote by $(\Gamma_{1}^-)_\lambda$ the set
 $\lambda\Gamma_{1}^-$. For any $y \in (\Gamma_{1}^-)_\lambda$, we define the function
  \[
x \in \ol \longmapsto  e^{i\tau_{\lambda, y} (x)} \eta\left(\frac{|x - y|}{\lambda}\right) w(x - y),
 \]
 where $\tau_{\lambda, y} (x) \doteq \sum\limits^{N}_{j = 1} A_\lambda^j(y)x^j$. {{By definition of $\eta$, this function belongs to $H^1_{A_\lambda}(\ol, \Gamma_{0 \lambda})$. From $(f_1)-(f_4)$, there exists $t_{\lambda, y}>0$  such that}}
\[
t_{\lambda, y}e^{i\tau_{\lambda, y}} \eta\left(\frac{| \cdot - y|}{\lambda} \right) w(\cdot - y) \in M_{\lambda}.
\]
{{Hence, $y \in (\Gamma_{1}^-)_\lambda$, and so we are able to  define the function}} $\Phi_\lambda : (\Gamma_{1}^-)_\lambda \to M_{\lambda}$ by
\begin{equation}\label{phi}
\Phi_\lambda (y)(x) = t_{\lambda, y} e^{i\tau_{\lambda, y} (x)} \eta\left(\frac{|x - y|}{\lambda}\right) w(x - y), \quad \forall x \in \ol.
\end{equation}

\begin{proposition} \label{phi}
Suppose that $f$ safisties $(f_1)$ and $(f_2)$. Then, the limit holds:
\[
\lim_{\lambda \to +\infty}\max_{y \in \Gamma_{1\lambda}^-}|I_{\lambda} (\Phi_\lambda (y))  - \ci| = 0.
\]
\end{proposition}

\dem\   Let $(\lambda_n)$ be any sequence such that $\lambda_n \to \infty$, as $n\to \infty$.  Since $(\Gamma_{1}^-)_{\lambda_n}$ is a compact set and $I_{\lambda_n}(\Phi_{\lambda_n}) \in C((\Gamma_{1}^-)_{\lambda_n})$, it suffices to prove that
\[
\lim_{n \to \infty} I_{\lambda_n} (\Phi_{\lambda_n} (y_n)) = c_\infty,
\]
for $y_n \in (\Gamma_{1}^-)_{\lambda_n}$   where the function
 $|I_{\lambda_n} (\Phi_{\lambda_n} (\cdot)) - \ci|$  attains its maximum {on $(\Gamma_{1}^-)_{\lambda_n}$}.

 By definition of $\grad_{A_{\lambda_n}}$,  for any $y \in (\Gamma_{1}^-)_{\lambda_n}$, we have
\[
|\grad_{A_{\lambda_n}} \Phi_{\lambda_n} (y)|^2 = \sum^N_{j = 1} |D^j_{A_{\lambda_n}} (\Phi_{\lambda_n} (y))|^2,
\]
where $D^j_{A_{\lambda_n}} (\Phi_{\lambda_n} (y) (x))$  are defined for $x \in \ol$ by
\begin{eqnarray*}
D^j_{A_{\lambda_n}} (\Phi_{\lambda_n} (y) (x))  & = & -i \partial_j \Phi_{\lambda_n} (y) (x) - A^j \left(\frac{x}{{\lambda_n}}\right)\Phi_{\lambda_n} (y) (x)\\
                                        & = & t_{{\lambda_n}, y} e^{i \tau_{{\lambda_n}, y (x)}} \left[ \eta\left(\frac{|x - y|}{{\lambda_n}}\right) w(x - y) \left( A^j\left(\frac{y}{{\lambda_n}}\right) - A^j\left(\frac{x}{{\lambda_n}}\right) \right)  \right.\\
                                        &   & \left. - i \partial_j \left( \eta \left(\frac{|x - y|}{{\lambda_n}}\right) w(x - y) \right) \right].
\end{eqnarray*}
Hence,
\begin{eqnarray*}
|\grad_{\al} \Phi_{\lambda_n} (y)(x)|^2 & = & t_{{\lambda_n}, y}^2 \left[ \eta^2 \left(\frac{|x - y|}{{\lambda_n}}\right) w^2(x - y) \left|A\left(\frac{y}{{\lambda_n}}\right)
                                 - A\left(\frac{x}{{\lambda_n}}\right) \right|^2  \right. \\
                                 &   & + \left. \left| \grad \left( \eta \left(\frac{|x - y|}{{\lambda_n}}\right) w(x - y) \right) \right|^2 \right].
\end{eqnarray*}
{Thereby}, for any $y \in (\Gamma_{1}^-)_{\lambda_n}$,
\begin{eqnarray*}
I_{{\lambda_n}} (\Phi_{\lambda_n} (y)) & = &\frac{t_{{\lambda_n}, y}^2}{2} \int_{\ol_n}\left\{ \eta^2\left(\frac{|x - y|}{{\lambda_n}} \right) w^2 (x - y) \left|A\left(\frac{y}{{\lambda_n}}\right) - A\left(\frac{x}{{\lambda_n}}\right)\right|^2 + \right.\\
                              &   & \left.+ \left|\grad \left(\eta\left( \frac{|x - y|}{{\lambda_n}}\right) w(x - y)\right) \right|^2 + \eta^2 \left( \frac{|x - y|}{{\lambda_n}}\right) w^2(x - y)\right\}dx - \\
                              &   & - \int_{\ol_n} F\left( t_{{\lambda_n}, y}^2 |\eta \left( \frac{|x - y|}{{\lambda_n}}\right) w(x - y)|^2\right)dx.
\end{eqnarray*}
Let $T_y$ be an orthogonal operator on $\rn$  which represents a rotation   such that the unitary normal vector to $T_y (\Omega_{\lambda_n} - \{y\})$ is $e_N=(0, \dots, 1)$. Set {{$\tilde \Omega_{{\lambda_n}, y} \doteq T_y (\Omega_{\lambda_n} - \{y\})$}}. After the change of variable $z=x-y$ and using that {{$\eta(\frac{| \cdot |}{\lambda_n}) $}} and $w$ are radially symmetric and  $T_y$ is a rotation, we find
\begin{eqnarray}
I_{{\lambda_n}} (\Phi_{\lambda_n} (y)) & = & \frac{t_{{\lambda_n}, y}^2}{2} \int_{\tilde \Omega_{{\lambda_n}, y}}\left[ \eta^2\left(\frac{|z|}{{\lambda_n}} \right) [w^2 (z) + |(\grad
                              w)(z)|^2] \right] \dd z \nonumber\\
                              &  - &  \frac{1}{2} \int_{\tilde\Omega_{{\lambda_n}, y}} F\left( t^2_{{\lambda_n}, y} \eta^2 \left( \frac{|z|}{{\lambda_n}} \right) w^2(z) \right) \dd z  \nonumber\\
                              &  +  & \label{1} \frac{t_{{\lambda_n}, y}^2}{2} \int_{\tilde \Omega_{{\lambda_n}, y}} \left|A\left(\frac{y}{{\lambda_n}}\right) - A\left(\frac{T^{-1}_y z + y}{{\lambda_n}}\right)\right|^2 \eta^2 \left( \frac{|z|}{{\lambda_n}} \right) w^2(z) \dd z  \\
                              & +  & \label{2}  \frac{t_{{\lambda_n}, y}^2}{2} \int_{\tilde \Omega_{{\lambda_n}, y}} \frac{1}{{\lambda_n}^2} \left| \eta' \left( \frac{|z|}{{\lambda_n}}\right)\right|^2 w^2(z) \dd z  \\
                              & +  & \label{3}  \frac{t_{{\lambda_n}, y}^2}{2} \int_{\tilde \Omega_{{\lambda_n}, y}} \frac{2}{{\lambda_n}} \left| \eta' \left( \frac{|z|}{{\lambda_n}} \right)\right| w(z) \eta \left( \frac{|z|}{{\lambda_n}} \right) |\grad w (z)| \dd z.
\end{eqnarray}
We claim that the respective integrals in \eqref{1}, \eqref{2} and \eqref{3} go to zero as ${n} \to +\infty$. Indeed, we first examine \eqref{1}. {Since $w \in L^{2}(\mathbb{R}^{N}),$} there exists $M > 0$ such that
\begin{equation}\label{parte1}
\int_{\tilde \Omega_{{\lambda_n}, y} \cap B^c_M(0)} \left|A\left(\frac{y}{{\lambda_n}}\right) - A\left(\frac{T^{-1}_y z + y}{{\lambda_n}}\right)\right|^2 \eta^2 \left( \frac{|z|}{{\lambda_n}} \right) w^2(z) \dd z < \epsilon.
\end{equation}
On the other hand, since $A$ in uniformly continuous on the compact set $\overline\Omega$, there exists  $\gamma > 0$ such that
\begin{equation}
\label{phi2}
|A(x + v) - A(x)| < \epsilon, \quad \forall |v| \le \gamma, \quad \forall x \in \overline\Omega.
\end{equation}
 Since $|T^{-1}_y z| \le M$ for all $z \in \tilde \Omega_{{\lambda_n}, y} \cap B_M(0)$, there exists  ${\lambda_n} > 0$ sufficiently large such that
$\left|{T^{-1}_y z}/{{\lambda_n}}\right| \le \gamma$, hence that, by   \eqref{phi2}, we have
\[
\left|A\left(\frac{y + T^{-1}_y z}{{\lambda_n}}\right) - A\left(\frac{y}{{\lambda_n}}\right)\right| < \epsilon, \quad \forall y \in (\Gamma_{1}^-)_{\lambda_n},
\]
for every ${\lambda_n} > 0$ sufficiently large. Thus, for every $ z \in B_M(0)$,
\[
\eta^2\left(\frac{|z|}{{\lambda_n}}\right) w^2 \left|A\left(\frac{y + T^{-1}_y z}{{\lambda_n}}\right) - A\left(\frac{y}{{\lambda_n}}\right)\right|^2\chi_{_{\tilde \Omega_{{\lambda_n}, y} \cap B_M(0)}} (z) \le {\epsilon^2}|w|^2_{\infty, \rn},
\]
and so
\begin{equation}
\label{phi1}
\lim_{{\lambda_n} \to \infty} \int_{\tilde \Omega_{{\lambda_n}, y} \cap B_M(0)} \left|A\left(\frac{y}{{\lambda_n}}\right) - A\left(\frac{T^{-1}_y z + y}{{\lambda_n}}\right)\right|^2 \eta^2 \left( \frac{|z|}{{\lambda_n}} \right) w^2(z) \dd z = 0.
\end{equation}
Combining \eqref{parte1} with \eqref{phi1}, gives that the integral in \eqref{1} goes to zero as ${\lambda_n} \to \infty$.
In order to analyze the integrals in \eqref{2}-\eqref{3}, take a constant $C > 0$ such that
\[
\chi_{\tilde\Omega_{{\lambda_n}, y}} \frac{1}{{\lambda_n}} \left( \frac{1}{2{\lambda_n}} \left| \eta' \left( \frac{|.|}{{\lambda_n}} \right) \right|^2 w^2 \right. + \left. \left| \eta' \left( \frac{|.|}{{\lambda_n}} \right) \right| \eta \left( \frac{|.|}{{\lambda_n}} \right) w |\grad w| \right)  \le C \left[ w^2 + w |\grad w| \right] \in L^1 (\rn)
\]
and
\begin{eqnarray*}
\chi_{\tilde\Omega_{{\lambda_n}, y}}(z) \frac{1}{{\lambda_n}} \left( \frac{1}{2{\lambda_n}} \left| \eta' \left( \frac{|z|}{{\lambda_n}} \right) \right|^2 w^2(z) + \left| \eta' \left( \frac{|z|}{{\lambda_n}} \right) \right| \eta \left( \frac{|z|}{{\lambda_n}} \right) w(z) |\grad w (z)| \right) \\
                                    \le \frac{C}{{\lambda_n}} {\left[|w(z)|^2 +  |w(z)| |\grad w(z)|\right]\to 0},
\end{eqnarray*}
almost everywhere $z \in \mathbb{R}^N$, as ${n} \to \infty$. By Lebesgue's dominated converge theorem, it follows that the integrals in \eqref{2} and \eqref{3} go to zero as  ${\lambda_n} \to +\infty$. Consequently,
\begin{eqnarray}
I_{{\lambda_n}} (\Phi_{\lambda_n} (y)) & = & \frac{t_{{\lambda_n}, y}^2}{2} \int_{\tilde \Omega_{{\lambda_n}, y}} \eta^2\left(\frac{|z|}{{\lambda_n}} \right) [|w|^2 + |\grad w|^2] \dd z  \nonumber\\
                              &   & - \frac{1}{2} \int_{\tilde\Omega_{{\lambda_n}, y}} F\left( t^2_{{\lambda_n}, y} \eta^2 \left( \frac{|z|}{{\lambda_n}} \right) |w|^2 \right) \dd z + t^2_{{\lambda_n}, y} o_{\lambda_n}(1), \label{phi3}
\end{eqnarray}
where $o_{\lambda_n}(1)$ denotes a quantity going to zero as ${n} \to \infty$. Taking $y=y_n$ and using  the notation,
 $\Omega_n = \Omega_{{\lambda_n}}, \tilde \Omega_n = \tilde \Omega_{{\lambda_n}, y_n}, t_n = t_{{\lambda_n}, y_n}$, we get
\begin{eqnarray}
I_{{\lambda_n}} (\Phi_{{\lambda_n}} (y_n)) & = & \frac{t_{n}^2}{2} \int_{\tilde \Omega_{n}} \eta^2\left(\frac{|z|}{{\lambda_n}} \right) [|(\grad w)(z)|^2 + w^2 (z)] \dd z - \nonumber\\
                                                &  & - \frac{1}{2} \int_{\tilde\Omega_{n}} F\left( t^2_{n} \eta^2 \left( \frac{|z|}{{\lambda_n}} \right) w^2(z) \right) \dd z + o_n(1) t^2_n \label{i1}
\end{eqnarray}
We claim that $t_n \to 1$, as $n\to \infty$. In fact, combining the definition of $t_n$ with the argument used in the study of the integrals \eqref{1}-\eqref{3}, yields
\begin{align}
o_n (1) + \int_{\tilde \Omega_{n}} \eta^2\left(\frac{|z|}{{\lambda_n}} \right) [|\grad w|^2 + &w^2] \dd z = \int_{\tilde\Omega_{n}} f\left( t^2_{n} \eta^2 \left( \frac{|z|}{{\lambda_n}} \right) w^2\right) \eta^2 \left( \frac{|z|}{{\lambda_n}} \right) w^2 \dd z .\label{tn}
\end{align}
To establish the boundedness of $(t_n)$, suppose by contradiction that there exists a subsequence  $t_{n_i} \to +\infty$. {{Using $(f_3)-(f_5)$, \eqref{CEEIU}, Fatou lemma,  $w >0$ in $\rn$  and  \eqref{tn}, we have
\begin{eqnarray*}
+\infty & > & \int_{\rn} (|\grad w|^2 + w^2)  =  \lim_{i \to \infty} \int_{\tilde \Omega_{n_i}} \eta^2 \left( \frac{|z|}{{\lambda_{n_i}}} \right) (|\grad w (z)|^2 + w^2 (z)) \dd z \\
       & \ge & \lim_{i \to \infty} \int_{{{B_{r_0} (r_0 e_N)}}} f \left( t^2_n \eta^2 \left( \frac{|z|}{{\lambda_{n_i}}} \right) w^2 (z) \right) \eta^2 \left( \frac{|z|}{{\lambda_{n_i}}} \right) w^2 (z) \dd z \\
       & = & \lim_{i \to \infty} \int_{{B_{r_0} (r_0 e_N)}} f \left( t^2_{n_i} w^2 (z) \right) w^2 (z) \dd z \\
       & = & +\infty,
\end{eqnarray*}}}
which is impossible. Hence, $(t_n)$ is a bounded sequence.  We can clearly assume that $t_n \to t_0$, as $n\to \infty$. To verify that $t_0 >0$, suppose by contradiction that $t_0 = 0$. By $(f_{1})-(f_{2})$ and Lebesgue's dominated convergence, we obtain
\begin{equation}\label{auxiliar}
\lim_{n \to \infty} \int_{\tilde \Omega_{n}} f \left( t^2_{n}\eta^2 \left( \frac{|z|}{{\lambda_{n}}} \right) w^2(z) \right) \eta^2 \left( \frac{|z|}{{\lambda_{n}}} \right) w^2(z) \dd z = 0.
\end{equation}
{On the other hand, from \eqref{auxiliar},  \eqref{tn} and \eqref{CEEIU},}  we have
\begin{align*}
0 & < \int_{\rn}( |\grad w|^2 + w^2) \dd z = \lim_{n\to \infty} \int_{\tilde \Omega_{n}} \eta^2 \left( \frac{|z|}{{\lambda_{n}}} \right) (|\grad w (z)|^2 + w^2 (z)) \dd z \\
  & = \lim_{n\to \infty} \int_{\tilde \Omega_n} f \left( t^2_{n} \eta^2 \left( \frac{|z|}{{\lambda_{n}}} \right) w^2 (z) \right) \eta^2 \left( \frac{|z|}{{\lambda_{n}}} \right) w^2 (z) \dd z = 0,
\end{align*}
{which is a contradiction.} Hence,  $t_n  \to t_0 >0$, as $n\to \infty$.  Now observe that
\begin{eqnarray*}
\int_{\rn} (|\grad w|^2 + w^2) dz& = & \lim_{n \to \infty} \int_{\tilde \Omega_n} \eta^2 \left( t_n^2\frac{|z|}{{\lambda_n}} \right) [|\grad w(z)|^2 + w^2(z)] \dd z \\
                             & = & \lim_{n \to \infty} \int_{\tilde \Omega_n} f\left( t_n^2 \eta^2 \left( \frac{|z|}{{\lambda_n}} \right) w^2(z) \right) \eta^2 \left( \frac{|z|}{{\lambda_n}} \right) w^2(z) \dd z\\
                             & = & \int_{\rn} f( t_0^2 w^2) w^2 \dd z.
\end{eqnarray*}
Using this, $(f_4)$ and  the properties on  $w$,  we conclude that $t_0 = 1$ . Therefore, the proposition follows from  \eqref{i1} and the Lebesgue's dominated convergence.
 \cqd

{Finally, we establish  a version of Lions's lemma \cite{lions}, whose proof proceeds along the same lines as in \cite[Lemma 2.1]{ZQWang} combined with interpolation of the $L^{p}$ spaces.
\begin{lemma} \label{lionsnosso}Let $l>0$, $2 \le s < 2^*$ and $\lambda_n \to +\infty$. Let $\{u_{n}\} \subset   H^1 (\Omega_{\lambda_n})$ be a  sequence such that
\[
\lim_{n \to \infty} \sup_{y \in \rn} \int_{B_l (y) \cap \Omega_{\lambda_n} } |u_{n}|^s \dd x= 0.
\]
Then
\[
\lim_{n \to \infty} \int_{ \Omega_{\lambda_n}} |u_{n}|^m \dd x= 0,
\]
for every $m \in (2, 2^*)$.
\end{lemma}}

\section{The behavior of the minimax levels}
Taking  $b_\lambda$ given by \eqref{blambda}, we have:
\begin{proposition} \label{bltoci}
$\displaystyle \lim_{\lambda \to \infty}b_{\lambda} = \ci.$
\end{proposition}
The proof of Proposition \ref{bltoci} is long and will be carried out in a series of steps. First, by definition of $\Phi_\lambda (y)$ and  Proposition \ref{phi},
\begin{equation}
\label{bltoci0}
b_{\lambda} \le I_{\lambda} (\Phi_{\lambda} (y)) = o_\lambda(1) + \ci.
\end{equation}
We now consider the auxiliary problems:
\begin{equation}
\left\{\begin{array}{rcll}
            -\Delta u + u & = & f(u^2) u  &\textrm{ in } \Omega_\lambda, \\
            \displaystyle\frac{\partial u}{\partial \nu} & = & 0  &\textrm{ on } \Gamma_{1 \lambda},\\
            u & = & 0 & \textrm{ on } \Gamma_{0 \lambda},
        \end{array}
\right. \label{2pl}
\end{equation}
and
\begin{equation}
\left\{\begin{array}{rcll}
            -\Delta u + u & = & f(u^2) u & \textrm{ in } \Omega_\lambda, \\
            \displaystyle\frac{\partial u}{\partial \nu} & = & 0 &\textrm{ on } \partial\Omega_{\lambda}.
        \end{array}
\right.  \label{2plb}
\end{equation}
We will denote by $H^1(\ol, \Gamma_{0 \lambda})$ be the Hilbert space
\[
H^1(\ol, \Gamma_{0 \lambda}) \doteq \left\{ u \in H^1(\ol); \textrm{ trace of } u = 0 \textrm{ on } \Gamma_{0 \lambda} \right\},
\]
endowed with the norm
\[
\| u \|_{\Omega_\lambda} = \left( \int_{\ol}(|\grad_\al u |^2 + |u|^2)  \dd x  \right)^{1/2}.
\]
Let $J_\lambda: H^1(\ol, \Gamma_{0 \lambda}) \to \mathbb{R}$ be the functional associated with \eqref{2pl} and given by
\[
J _\lambda (u) = \frac{1}{2} \int_{\ol} (|\grad u|^2 + u^2)\dd x - \frac{1}{2}\int_{\ol} F(u^2)\dd x, \quad \forall u \in H^1(\ol, \Gamma_{0 \lambda}).
\]
We define the functional $\overline J _\lambda: H^1(\ol) \to \re$ associated with \eqref{2plb} by
\[
\overline J _\lambda (u) = \frac{1}{2} \int_{\ol} (|\grad u|^2 + u^2)\dd x - \frac{1}{2}\int_{\ol} F(u^2)\dd x, \quad \forall u \in H^1(\ol),
\]
with corresponding Nehari manifold  and mountain pass level given by
\[
\overline N_\lambda\doteq \{ u \in H^1(\ol) \backslash \{0\};\ \overline J'_{\lambda}(u)u = 0\}
\quad \mbox{and}\quad
\overline c_{\lambda}  \doteq \inf\limits_{\overline N_\lambda} \overline J_{\lambda}.
\]
We will also denote by  $c _\lambda$ the  mountain pass level associated with the problem \eqref{2pl}. By the definition of these levels and from \eqref{diamagnetic}, we find
\begin{equation}
\label{bltoci00}
b_{\lambda} \ge c_{\lambda} \ge \overline c_{\lambda} > 0.
\end{equation}
From \eqref{bltoci0}-\eqref{bltoci00}, we deduce that it suffices to show that
\begin{equation}
\label{bltoci1}
\lim_{\lambda \to \infty}\overline c_{\lambda} = \ci.
\end{equation}
In order to prove \eqref{bltoci1}, we begin by observing  that the  mountain pass theorem combined with a similar argument employed in the proof of Proposition \ref{ps} implies that  there is a solution $u_\lambda \in H^1(\ol)$ of \eqref{2plb} satisfying
\begin{equation}
\label{bltoci10}
\overline J_{\lambda} (u_\lambda) = \overline c_{\lambda} = \inf_{\overline N_\lambda} \overline J_{\lambda}, \quad \overline J'_{\lambda}(u_\lambda) = 0,
\end{equation}
for every $\lambda >0$.  Combining \eqref{bltoci00} with \eqref{bltoci10}, gives that $\sup_{\lambda >0}\overline J_\lambda (u_\lambda) < \infty$ and  $\overline J'_\lambda (u_\lambda) u_\lambda = 0$ for all $\lambda > 0$. By $(f_3)$,
\begin{equation}
\label{bltoci2}
\sup_{\lambda >0}\|u_\lambda\|_\ol < \infty
\end{equation}
(where $\|\cdot \|_\ol$ denotes the norm of $H^1(\ol)$). Exploiting similar argument used in the proof of
Proposition \ref{limnehari}, we may assume that
\begin{equation}
\label{bltoci3}
\| u_\lambda \|^2_{\ol} \ge \delta_0 \quad \textrm{and} \quad \overline J_{\lambda} (u_\lambda) = \overline c _\lambda \ge \delta_0, \quad \forall\, \lambda > 0, 
\end{equation}
for some constant $\delta_0 > 0$ independent of $\lambda$.  From  \eqref{bltoci3} and Lemma \ref{lionsnosso}, there exist $(y_\lambda)_\lambda \subset \rn, l > 0$ and $\gamma > 0$ such that
\begin{equation}
\label{bltoci4}
\liminf_{\lambda \to \infty} \int_{\ol \cap B_l(y_\lambda)} |u_\lambda|^2\dd x \ge \gamma > 0.
\end{equation}
Moreover, by increasing $l$ if necessary, we may assume that  $y_\lambda \in \ol$ for every $\lambda > 0$, because \eqref{bltoci4} yileds $\ol \cap B_l(y_\lambda) \neq \emptyset$, for every $\lambda$.

\begin{lemma}\label{fronteira} There exists a constant $C>0$  such that  $\textrm {dist}(y_\lambda, \partial \ol) \le C$, for every $\lambda >0$. \end{lemma}

\dem\ Suppose the lemma were false. {Then, we could find a sequence} $(\lambda_n)$ such that $\lambda_n \to \infty$ and  $\textrm {dist}(y_{\lambda_n}, \partial \Omega_{\lambda_n}) \to  \infty$, as $n\to \infty$. Let $R > l$ be an arbitrary number. For $n$ sufficiently large, we have $B_{2R}(y_{\lambda_n}) \subset \Omega_{\lambda_n}$. Define
\[
w_{{\lambda_n}, R} (x) \doteq \eta\left(\frac{|x|}{R}\right) u_{\lambda_n} (x + y_{\lambda_n}), \quad \forall\, x \in \Omega_{\lambda_n}-\{y_{\lambda_n}\},
\]
where $\eta \in C^{\infty} (\re)$ is such that $\eta = 1, \textrm{ on } [0,1]$,  $\eta = 0, \textrm{ on } (2,+\infty)$, $0 \leq \eta  \leq 1$ and $\eta' \in L^{\infty}(\re)$. Hence,
$\supp \, w_{{\lambda_n}, R}\subset B_{2R}(0)$. We can assume that $w_{{\lambda_n}, R} \in H^1(\rn)$ and also $\sup_{n}\|w_{{\lambda_n}, R}\| \leq C$, for some constant $C>0$ independent $R$. Observing that
\[
\int_{B_l(0)}|w_{{\lambda_n}, R}|^2 \dd x = \int_{B_l(0)}|u_{{\lambda_n}} (x + y_{\lambda_n})|^2 \dd x = \int_{B_l(y_{\lambda_n})}|u_{{\lambda_n}}|^2 \dd x \ge \gamma > 0,
\]
we get a nontrivial function $w_R \in H^1(\rn)$ such that
\[
\left\{ \begin{array}{l}
w_{{\lambda_n}, R} \rightharpoonup w_R,   \textrm{ weakly in } H^1(\rn), \textrm{ as } n\to \infty, \\

w_{{\lambda_n}, R} \to w_R,    \textrm{ strongly in } L^p_{loc}(\rn),  p \in \left.[1, 2^*\right.),  \textrm{ as } n\to \infty,    \\

\int_{B_l (0)}|w_R|^2 \ge \gamma > 0.
\end{array}
\right.
\]
Let $\|\cdot \|$ denote the norm in of $H^1(\rn)$. {{Since $\|w_R\| \le \liminf_{{n} \to \infty} \|w_{{\lambda_n}, R}\|$}},
the family $(w_R)_R \subset H^1(\rn)$ is bounded. Hence, there exists $v \in H^1(\rn)$ such that
\[
\left\{ \begin{array}{l}
w_R \rightharpoonup v, \textrm{ weakly in } H^1(\rn), \textrm{ as } R \to \infty,\\

w_{R} \rightarrow v, \textrm{ strongly in } L^p_{loc}(\rn), p \in \left.[1, 2^*\right.), \textrm{ as } R\to \infty,\\

\int_{B_l (0)}|v|^2 \ge \gamma > 0.
\end{array}
\right.
\]
In particular, $v \not\equiv 0$. We assert that  $v$ is a solution of \eqref{2prn}. In fact, given $\phi \in \cinfsc (\rn)$, we take $t > 0$ such that  $\supp \phi \subset B_t(0)$ and $B_t(y_{\lambda_n}) \subset \Omega_{\lambda_n}$ for $n$ sufficiently large. As $u_{\lambda_n}$ is a weak solution of \eqref{2plb} for $\lambda = \lambda_n$, we have
\begin{eqnarray*}
{\int_{B_t(0)} [\grad u_{\lambda_n} (x + y_{\lambda_n}) \grad \phi + u_{\lambda_n} (x + y_{\lambda_n}) \phi ]} & = & \int_{\Omega_{\lambda_n}} [\grad u_{\lambda_n} (x + y_{\lambda_n}) \grad \phi + u_{\lambda_n} (x + y_{\lambda_n}) \phi]\\
                                                                                       & = & \int_{\Omega_{\lambda_n}} f(u^2_{\lambda_n} (x + y_{\lambda_n})) u_{\lambda_n} (x + y_{\lambda_n}) \phi\\
                                                                                       & = & \int_{B_t(0)} f(u^2_{\lambda_n} (x + y_{\lambda_n})) u_{\lambda_n} (x + y_{\lambda_n}) \phi.
\end{eqnarray*}
For $n$ sufficiently large and $R > t$, we obtain
\[
\int_{B_t(0)}[ \grad w_{{\lambda_n}, R} \grad \phi + w_{{\lambda_n}, R} \phi]\dd x = \int_{B_t(0)} f(w^2_{{\lambda_n}, R}) w_{{\lambda_n}, R} \phi\dd x.
\]
Taking $n \to \infty$, we have
\[
\int_{B_t(0)}[ \grad w_{R} \grad \phi + w_{R} \phi ]\dd x= \int_{B_t(0)} f(w^2_{R}) w_{ R} \phi\dd x.
\]
Using that $\supp \phi \subset B_t(0)$ and  $R > t$, we find after taking $R \to \infty$
\[
\int_{\rn} [ \grad v \grad \phi + v \phi]  = \int_{B_t(0)} [\grad v \grad \phi + v \phi] = \int_{B_t(0)} f(v^2) v \phi = \int_{\rn} f(v^2) v \phi.
\]
Since $\phi \in \cinfsc (\rn)$ is arbitrary, we conclude that $v$ is a nontrivial solution of  \eqref{2prn}. Given $M > R$, we take $n$ sufficiently large such that $B_M(y_{\lambda_n}) \subset \Omega_{\lambda_n}$. By  \eqref{bltoci0}-\eqref{bltoci00},
\begin{eqnarray*}
o_{\lambda_n}(1) + \ci & \ge & \overline c_{{\lambda_n}} = \overline J_{{\lambda_n}} (u_{\lambda_n}) - \frac{1}{2}\overline J_{{\lambda_n}}'(u_{\lambda_n}) u_{\lambda_n} \\
             & = & \frac{1}{2} \int_{\Omega_{\lambda_n}} [f(u_{\lambda_n}^2)u_{\lambda_n}^2 - F(u^2_{\lambda_n})]\dd x  \\
             & \ge & \frac{1}{2} \int_{B_M(y_{\lambda_n})} [f(u_{\lambda_n}^2)u_{\lambda_n}^2 - F(u^2_{\lambda_n})]\dd x\\
             & = & \frac{1}{2}\int_{B_M(0)} [f(w_{{\lambda_n}, R}^2)w_{{\lambda_n}, R}^2 - F(w_{{\lambda_n}, R}^2)]\dd x.
\end{eqnarray*}
By Fatou's lemma and \eqref{crncinf}, we obtain, after taking  ${n} \to \infty$,  $R \to \infty$ and $M \to \infty$,
\[
\ci \ge \frac{1}{2} \int_{\rn} [f(v^2) v^2 - F(v^2)]\dd x = J_{\rn}(v) \ge c_{\rn} = 2\ci,
\]
 which is a contradiction. Lemma \ref{fronteira} is proved.  \cqd

From Lemma \ref{fronteira}, by increasing  $l$ if necessary, we may assume that $y_\lambda \in \partial \ol$ in \eqref{bltoci4}.  Let $T_{y_\lambda}$ be an orthogonal operator on $\rn$  which represents a rotation  such that the inward unitary normal vector to
{{$\tilde\Omega_\lambda \doteq T_{y_\lambda} \left(\Omega_\lambda -\{ y_\lambda\} \right) $}} is $e_N=(0, \cdots, 1)$.   We define
\[
v_\lambda (x) = u_\lambda (T^{-1}_{y_\lambda} x + y_\lambda ), \quad \forall x \in \tilde\Omega_\lambda.
\]
In the following, we gather the properties satisfied by $v_\lambda$:
\begin{itemize}
\item[$(a)$] Since $\|v_\lambda\|_{\tilde\Omega_\lambda} = \|u_\lambda\|_{\ol}$,  \eqref{bltoci2} shows that  $\sup_{\lambda>0}\|v_\lambda\|_{\tilde\Omega_\lambda} < \infty$;
\item[$(b)$] $\displaystyle{\int_{\tilde\Omega_\lambda} F(v_\lambda^2) \dd x = \int_{\ol} F(u^2_\lambda)} \dd x$;
\item[$(c)$] Since $u_\lambda$ is a solution of \eqref{2plb}, $v_\lambda$ is a solution of
\begin{equation}
\left\{\begin{array}{rcl}
            -\Delta u + u & = & f(u^2) u \textrm{ in } \tilde\Omega_\lambda , \\
            \displaystyle\frac{\partial u}{\partial \nu} & = & 0 \textrm{ on } \partial\tilde\Omega _\lambda;
        \end{array}
\right.  \label{2potlb}
\end{equation}
\item[$(d)$] $J_{\tilde\Omega_\lambda}(v_\lambda) = \overline c_{\tilde\Omega_\lambda} = \overline c_{\lambda}$, where $J_{\tilde\Omega_\lambda}$ is the functional associated with \eqref{2potlb} and  $\overline c_{\tilde \Omega_\lambda}$ is the corresponding mountain pass level;
\item[$(e)$] From \eqref{bltoci4},
\[
\liminf\limits_{\lambda \to \infty} \displaystyle\int_{B_l(0) \cap \tilde\Omega_\lambda} |v_\lambda|^2 \ge \gamma.
\]
\end{itemize}
Given  $\rho > h > 0$, we define
\[
D_{\rho, h} \doteq \left\{ (x_1, \dots , x_N)  \in \mathbb{R}^N ;\ x_N > h \right\} \cap B_{\rho} (0).
\]
From \eqref{CEEIU}, $\chi_{\tilde\Omega_\lambda} \to \chi_{{\mathbb{R}_{+}^N}}$ almost everywhere in $\rn$, as $\lambda \to \infty.$ Hence,  $D_{\rho, h} \subset \tilde\Omega_\lambda$ for every $\lambda$ sufficiently large. Thus, $v_\lambda \in H^1(D_{\rho, h})$ for every $\lambda$ sufficiently large. By $(a)$, we may assume that there exists $v_{\rho, h} \in H^1(D_{\rho, h})$ such that
\[
\left\{
\begin{array}{ll}
v_{\lambda} \rightharpoonup v_{\rho, h} & \textrm{ weakly in } H^1(D_{\rho, h}), \textrm{ as } \lambda \to \infty, \\

v_{\lambda} \to  v_{\rho, h} & \textrm{ strongly in  } L^p(D_{\rho, h}), p \in \left.[1, 2^*\right.), \textrm{ as } \lambda \to \infty,  \\

v_{\lambda} (x) \to v_{\rho, h}(x) &  \textrm{ a.e. in } D_{\rho, h}, \textrm{ as } \lambda \to \infty.
\end{array}
\right.
\]
Using $(a)$ one more time and the Banach-Steinhaus theorem, we find a constant $K>0$ such that
\[
\|v_{\rho, h}\| _{D_{\rho, h}} \le K, \quad \forall\, \rho, h > 0
\]
(where $\|\cdot \|_{D_{\rho, h}}$ denotes the norm of $H^1(D_{\rho, h})$).
Let $\rho_n \to \infty$ and $h_n \to 0$ be monotone sequences. Thus,   $$D_n \doteq D_{\rho_n, h_n} \subset D_{\rho_{n+1}, h_{n+1}} \doteq D_{n+1},\quad \forall n \ge 1.
$$
{ This  allows us to apply a diagonal type argument to obtain a  bounded subsequence $(v_{k})$ in  $H^1(\rn_+)$  and a function $v \in H^1(\rn_+)$ such that
\begin{eqnarray} \label{bltoci40}
& v_{k} \rightharpoonup  v \textrm{ weakly in  } {{H^1(\rn_+)}}, \textrm{ as } k \to \infty, \nonumber \\
& v_{k} \to  v \textrm{ strongly in  } L^p_{loc}(\rn_+), \forall p \in \left.[1, 2^*\right.), \textrm{ as } k \to \infty,\\
& v_{k}(x) \to v(x) \textrm{ a.e. in } \rn_+, \textrm{ as } k \to \infty. \nonumber
\end{eqnarray}}
\begin{lemma}\label{weaklimit}
The function $v$ is a nontrivial weak solution of \eqref{2pinf}.
\end{lemma}

\dem\  We first show that  $v\not \equiv 0$.  In fact,
from $(e)$,
\begin{equation}
\label{bltoci5}
\liminf_{k \to \infty} \int_{B_l(0) \cap \tilde\Omega_{k}} v^2_{k} \ge \gamma > 0.
\end{equation}
 Given $t \in (0, l)$, define $A_t = \left\{ x \in B_l(0) \cap \tilde\Omega_{k}; 0 \le x^N \le t \right\}$ and  $\Lambda_{k} = \left(B_l(0) \cap \tilde\Omega_{k}\right) \backslash A_t$. Thus,
\[
\int_{B_l(0) \cap \tilde\Omega_{k}} v^2_{k} = \left(\int_{A_t} + \int_{\Lambda_{k}}\right) v^2_{k}.
\]
 As $\sup_{k}\|v_{k}\|_{D_k} < \infty$, using H\"older's inequality and the Sobolev embedding theorem, we get
\[
\int_{A_t} v_{k}^2 \le \left(\int_{A_t} v_{k}^{2^*}\right)^{\frac{2}{2^*}} \left(\int_{A_t} 1 \right)^{\frac{2}{N}} \le \overline K |A_t|^{\frac{2}{N}},
\]
for some constant $\overline K > 0$. Now choose a $t \in (0, l)$ such that
\[
\int_{A_t} v_{k}^2 \le \left(\int_{A_t} v_{k}^{2^*}\right)^{\frac{2}{2^*}} \left(\int_{A_t} 1 \right)^{\frac{2}{N}} \le \overline K |A_t|^{\frac{2}{N}} < \frac{\gamma}{4}.
\]
Consequently, from \eqref{bltoci5}, for all sufficiently large  $k$, we have
\[
\frac{\gamma}{2} \le \int_{B_l(0) \cap \tilde\Omega_{k}} v^2_{k} \le \frac{\gamma}{4} + \int_{\Lambda_{k}} v^2_{k} \le \frac{\gamma}{4} + \int_{D} v^2_{k},
\]
for every compact set $D \subset \rn$ with $\Lambda_{k} \subset D \subset D_{k}$. Hence, for all sufficiently large  $k$,
\[
\int_{D} v^2_k \ge \frac{\gamma}{4}
\]
and consequently
\[
\int_{D} v^2 = \lim_{k \to \infty} \int_{D} v^2_{k} \ge \frac{\gamma}{4} > 0,
\]
which implies $v \not \equiv 0$.  In order to prove that $v$ is a weak solution of \eqref{2pinf}, we first show that
$\grad v_k  \to \grad v$, strongly in $(L^2(K))^N$, for any compact set $K \subset \rn_+$.  Effectively, let    $K \subset \rn_ +$ be a compact set.  Taking
$\psi \in \cinfsc(\rn_+)$ such that $\psi \equiv 1$, on $K$, and $0\leq \psi \leq 1$, we have  $\supp\, \psi \subset \tilde\Omega_k$, for every  $k$  sufficiently large. As
$v_k \psi$, $v\psi \in H^1(\tilde\Omega_k)$ and  $v_k$ is a weak solution of  \eqref{2potlb}, we have
\begin{eqnarray}
\overline J '_{\tilde\Omega_k} (v_k) (v_k \psi) = \int_{\tilde\Omega_k}[ |\grad v_k |^2\psi + v_k \grad v_k \grad\psi + v_k^2 \psi] - \int_{\tilde \Omega_k} f(v_k^2)v_k^2 \psi &=& 0 \label{aux1}\\
\overline J '_{\tilde\Omega_k} (v_k) (v \psi) = \int_{\tilde\Omega_k}[ \psi \grad v_k \grad v + v \grad v_k \grad\psi + v_k v \psi] - \int_{\tilde \Omega_k} f(v_k^2)v_k v \psi &=& 0 \label{aux2}
\end{eqnarray}
where $\overline J_{\tilde \Omega_k}: H^1(\tilde \Omega_k) \to \re$ is the functional associated with \eqref{2potlb}. Combining \eqref{aux1}-\eqref{aux2}, we obtain
\begin{eqnarray*}
 \int_{K} |\grad v_k - \grad v|^2 &\le& \int_{\rn} \psi \left[ |\grad v_k|^2 - 2\grad v_k \grad v + |\grad v|^2\right] \\
   & = & \int_{\rn_+} [\psi |\grad v_k|^2 - \psi \grad v_k \grad v + \psi \grad v \grad( v - v_k) ]\\
   & = & \int_{\rn_+} [ f(v_k^2) v_k ^2 \psi - v_k \grad v_k \grad \psi - v_k ^2 \psi] + \int_{\rn_+} [v \grad v_k \grad \psi + v_k v \psi ]\\
    &   &  - \int_{\rn_+}f(v_k^2) v_k v \psi
  + \int_{\rn_+}\psi \grad v \grad( v -  v_k)\\
   & = & \int_{\rn_+} {[ f(v_k^2) \psi v_k (v_k - v) - (v_k - v) \grad v_k \grad \psi - v_k \psi (v_k - v) ]}\\
   &   &  +  \int_{\rn_+}\psi \grad v \grad( v_k -  v).
\end{eqnarray*}
This and the fact that $(v_{k})$ is bounded in $L^2(\rn_+)$ combined with $(f_{6})$,  \eqref{bltoci40} and  H\"older's inequality show that
\begin{eqnarray*}
\int_{K} |\grad v_k - \grad v|^2 \le o_k(1),\quad \mbox{as $k\to \infty$},
\end{eqnarray*}
that is $ \grad v_k \to \grad v$,  strongly in  $(L^2(K))^N$, as desired. As a  consequence,
\begin{equation}
\label{bltoci50}
\grad v_{k} (x) \to \grad v(x),\quad \textrm{ for  almost every  }  x\in \rn.
\end{equation}
In order to conclude the proof of Lemma \ref{weaklimit}, it remains to prove that
\begin{equation}\label{aux3}
\int_{\rn_+} \left[\grad v \grad \phi + v \phi\right] - \int_{\rn_+} f(v^2)v\phi = 0,\quad \forall\, \phi \in H^1(\rn_+).
\end{equation}
Since the set of restrictions of the functions of  $\cinfsc(\rn)$ to $\rn_+$ is a dense subspace of $H^1(\rn_+)$ (see \cite[Corollaire IX.8]{brezis}), it suffices to show that relation \eqref{aux3} holds for every  $\phi \in \cinfsc(\rn)$. Given $\phi \in \cinfsc(\rn)$, {let $t > 0$ be such that $B_t(0) \supset \supp\, \phi$}. From \eqref{CEEIU}, $\chi_{\tilde \Omega_k \cap B_t(0)}  \to \chi_{B^+_t}$ almost everywhere in $\rn$, as $k \to \infty$,  where
 $B^+_t \doteq B_t(0) \cap \rn_+$ and where $\chi_{B^+_t}$ is the characteristic function related to the set $B^+_t$.  This and  \eqref{bltoci50} imply that
$\chi_{\tilde \Omega_k \cap B_t(0)} \grad v_k {\to} \chi_{B^+_t} \grad v$, almost everywhere in $\rn$, as $k \to \infty$. Furthermore, $(\chi_{\tilde \Omega_k \cap B_t(0)} \grad v_k)_k$ is bounded in \ $(L^2(\rn_+))^N$. Hence,
$\chi_{\tilde \Omega_k \cap B_t(0)} \grad v_k {\tofraco} \chi_{B^+_t} \grad v$  weakly in $(L^2(\rn_+))^N$,  as $k \to \infty$, and so
\begin{equation}
\label{bltoci51}
\lim_{k \to \infty}\int_{\tilde \Omega_k} \grad v_k \grad \phi = \lim_{k \to \infty}\int_{\rn_+} \chi_{\tilde \Omega_k \cap B_t(0)} \grad v_k \grad \phi =  \int_{\rn_+} \chi_{B^+_t} \grad v \grad \phi = \int_{\rn_+} \grad v \grad \phi.
\end{equation}
{Since $(v_k)$ is bounded in $H^1(\rn_+)$, by $(f_{6})$, there exists $M_1>0$ such that
\begin{equation}\label{M1}
\int_{B_t(0)} |f(v_k^2)v_k|^{q/(q-1)} \leq M_1.
\end{equation}
Given $\eta >0$, from \eqref{bltoci40} and Egoroff's theorem, there exists $E\subset B_t(0)$ such that $|E| < \eta$ and $v_k(x) \to v(x)$ uniformly on $ B_t(0) \setminus E$. Using H\"{o}lder's inequality, \eqref{M1} and $(f_{6})$, we get $M_2>0$ such that
\[
\left| \int_{B_t(0)}( f(v_k^2)v_k  - f(v^2)v) \phi\right| \leq
\int_{B_t(0) \setminus E}| f(v_k^2)v_k  - f(v^2)v| |\phi| + M_2\eta^q.
\]
As $\eta >0$ can be chosen arbitrarily small, $f(v_k^2)v_k \to f(v^2)v$ uniformly on $B_t(0) \setminus E$ and $\supp\, \phi \subset B_t(0)$,  we obtain }
\begin{equation}
\label{bltoci53}
\lim_{k \to \infty}\int_{\tilde \Omega_k} f(v_k^2)v_k \phi  =  \int_{\rn_+} f(v^2)v \phi.
\end{equation}
Using \eqref{bltoci40}, similar arguments to those above show that
\begin{equation} \label{bltoci52}
\lim_{k \to \infty} \int_{\tilde \Omega_k} v_k \phi = \int_{\rn_+} v \phi.
\end{equation}
Combing  \eqref{bltoci51} - \eqref{bltoci52} with the fact that  $v_k$ satisfies  \eqref{2potlb}, yields
\begin{eqnarray*}
0 = \lim_{ k \to \infty} \int_{\tilde \Omega_k} (\grad v_k \grad \phi + v_k \phi - f(v_k ^2) v_k \phi) = \int_{\rn_+} (\grad v \grad \phi + v \phi - f(v^2) v \phi),
\end{eqnarray*}
for every $\phi \in \cinfsc(\rn)$, and the proof is complete. \cqd

In the following, we conclude the proof of Propostion \ref{bltoci}. From  \eqref{bltoci0} and \eqref{bltoci00},
\begin{eqnarray*}
\ci + o_k (1) & \ge & \overline c_{\Omega_k} = \overline c_{\tilde \Omega_k} \\
               & = &  \overline J_{\tilde \Omega_k}(v_k) = \overline J_{\tilde \Omega_k}(v_k) - \frac{1}{2} \overline J'_{\tilde \Omega_k}(v_k)v_k \\
               & = & \frac{1}{2} \int_{\tilde\Omega_k}[ f(v_k^2) v_k^2 - F(v_k^2)].
\end{eqnarray*}
Using Fatou's lemma and  \eqref{bltoci40}, we have
\begin{eqnarray*}
\ci & \ge & \limsup_{k \to \infty} \overline c_{\tilde \Omega_k}
      \ge  \liminf_{k \to \infty} \overline c_{\tilde \Omega_k}
     =  \liminf_{k \to \infty} \frac{1}{2} \int_{\tilde\Omega_k} {[f(v_k^2) v_k^2 - F(v_k^2)]} \\
    & \ge & \frac{1}{2}\int_{\rn_+} {[f(v^2) v^2 - F(v^2)]}
     =  J_{\infty} (v) \ge  \ci.
\end{eqnarray*}
Consequently,
\[
\lim_{\lambda \to \infty} \overline c_{\Omega_\lambda} = \ci,
\]
that is, \eqref{bltoci1} holds, and the proof of Proposition \ref{bltoci} is complete. \cqd

\section{The barycenter map}
This section is devoted to establish a key relation between some subsets of $\mathbb{R}^N$ and $M_\lambda$. {For $q \in (2,2^{*})$ given by $(f_5)$ and $\lambda > 0$, define the barycenter map  $\beta_\lambda: M_{\lambda} \to \rn$  by}
\[
\beta_\lambda (u) = \frac{\displaystyle\int_{\Omega_\lambda} x |u|^q \dd x}{\displaystyle\int_{\Omega_\lambda} |u|^q \dd x}.
\]
\begin{proposition} \label{baricentro}
Let $(\Gamma_{1}^+)_\lambda$ be the expanding set $\lambda\Gamma_{1}^+$. Then,  there exist  $\epsilon^* > 0$ and  $\lambda_1 > 0$ such that $\beta_\lambda (u) \in (\Gamma_{1}^+)_\lambda$, provided that   $\lambda > \lambda_1$,   $u \in M_{\lambda}$ and $I_{\lambda} (u) \le b_\lambda^*$, where $b_\lambda^* = b_{\lambda} + \epsilon^*$.
\end{proposition}

\dem \ It suffices to show that if $(\epsilon_n)$ and $(\lambda_n)$ are arbitrary sequences, with $\epsilon_n \to 0$ and $\lambda_n  \to \infty$, and  if $u_n \in M_{\lambda_n}$ is a sequence such that
\begin{equation}
\label{baricentro1}
b_{\lambda_n} \le I_{\lambda_n} (u_n) \le b_{\lambda_n} + \epsilon_n ,
\end{equation}
then
\begin{equation}
\label{baricentro2}
\textrm{dist}(\beta_{\lambda_n}(u_n), \Gamma_{1\lambda_n }) \le \lambda_n r,
\end{equation}
for every $n$ sufficiently large.   In fact, by \eqref{baricentro1} and  Proposition \ref{bltoci},
\begin{equation}\label{baricentro22}
I_{\lambda_n} (u_n) \to \ci,\quad \mbox{as } n\to \infty.
\end{equation}
Using that $u_n \in M_{\lambda_n}$ and \eqref{diamagnetic}, there exists   $t_n > 0$ such that
\begin{equation}
\label{baricentro3}
o_n(1) + \ci = I_{\lambda_n} (u_n) \ge \max_{t \ge 0} I_{\lambda_n} (t u_n) \ge \max_{t \ge 0} J_{\lambda_{n}} (t |u_n|) = J_{\lambda_{n}} (t_n |u_n|) \ge c_{\lambda_{n}},
\end{equation}
where $J_{\lambda_{n}}$ is the functional associated with \eqref{2pl} with $\lambda=\lambda_n$, $c_{\lambda_{n}}$ and $N_{\lambda_{n}}$ are the corresponding  mountain pass level and the Nehari manifold.   Combining \eqref{bltoci00}-\eqref{bltoci1} with \eqref{baricentro3} and Proposition
 \ref{bltoci}, we have
\begin{equation}\label{baricentro33}
\lim_{n \to \infty} c_{\lambda_{n}} = \lim_{n \to \infty} J_{\lambda_{n}} (t_n |u_n|) = \ci.
\end{equation}
Set {$\overline\epsilon_n = b_{\lambda_n} - c_{\lambda_n}$}. By Proposition  \ref{bltoci} and \eqref{baricentro33}, $\overline\epsilon_n \to 0$, as $n\to \infty$. Thus,
\[
 \overline\epsilon_n + c_{\lambda_n}  \ge J_{\lambda_n} (t_n |u_n|) \ge c_{\lambda_n}.
\]
Applying Ekeland variational principle {\cite[Corollary 3.4]{Ekeland}}, for every  $\nat$,  there exists  $v_n \in N_{\lambda_n}$ such that
\begin{equation}
\label{baricentro4}
\|t_n |u_n| - v_n\|_{\Omega_n} \le 2 \sqrt{\overline\epsilon_n}, \quad c_{\lambda_n} \le J_{\lambda_n} (v_n) \le c_{\lambda_n} + 2 \overline \epsilon_n
\end{equation}
and
\[
\| \left( J_{\lambda_n}\big|_{{N_{\lambda_n}}} \right)'(v_n)\|_{(H^1(\Omega_{\lambda_n}, \Gamma_{0 \lambda_n}))'} \le 8 \sqrt{\overline \epsilon_n}.
\]
As in the proof of Proposition \ref{psr}, we find that  $v_n \in H^1(\Omega_{\lambda_n}, \Gamma_{0 \lambda_n})$ satisfies
\begin{equation}
\label{baricentro5}
J_{\lambda_n}(v_n) \to \ci, \quad J'_{\lambda_n}(v_n) \to 0.
\end{equation}
From \eqref{baricentro5} and $(f_3)$, the sequence $(\| v_n \|_{\Omega_{\lambda_n}})_n$ is bounded. Consequently, from Lemma
\ref{lionsnosso} and \eqref{baricentro5}, there exist  $l>0$, $\gamma > 0$ and  $y_n \in \rn$ such that
\[
\liminf_{n \to \infty}\int_{B_l(y_n) \cap \Omega_{\lambda_n}} |v_n|^2 \ge \gamma > 0.
\]
Proceeding as in the proof of Lemma \ref{fronteira}, with \eqref{baricentro5}  replacing \eqref{bltoci10}, we get a positive constant  $C>0$ such that $\textrm{dist} (y_n, \partial\Omega_{\lambda_n}) \le C$. Thus,  by increasing  $l$ if necessary, we may assume that   $y_n \in \partial\Omega_{\lambda_n}$.  Following the same argument in the proof of Proposition \ref{bltoci}, we define
\[
\tilde v_n (x) = v_n (T^{-1}_{y_n} x + y_n), \quad \forall x \in \tilde\Omega_n \doteq T_{y_n}(\Omega_{\lambda_ n} - \{y_n\}), \quad \forall \nat,
\]
to obtain a subsequence of  $\tilde v_n \in H^1(\tilde \Omega_n, \tilde\Gamma_{0 \lambda_n})$ (still denoted by $\tilde v_n$)  and a function $v \in H^1(\rn_+)$  such that
\begin{eqnarray}
& \tilde v_n \tofraco v \textrm{ weakly in } {{ H^1(\rn_+)}}, \, \tilde v_n \to v \textrm{ strongly in } L^p_{loc} (\rn_+), \forall p \in [\left.1, 2^*)\right. \label{baricentro7}\\
& \tilde v_n (x) \to v(x), \, \grad \tilde v_n (x) \to \grad v(x) \, \textrm{ almost every } x\in \rn_+. \label{baricentro8}
\end{eqnarray}

\smallskip
\noindent \textbf{Claim I.} \textit{There exits a constant  $C > 0$ such that $\textrm{ dist} (y_n, \Gamma_{1 \lambda_n}) \le C$.}

\smallskip
In fact, suppose Claim I were false. Then we could find subsequences (not renamed) such that
\begin{equation}
\label{baricentro6}
\alpha_n \doteq \textrm{dist}(y_n, \Gamma_{1 \lambda_n}) \to \infty ,\quad \mbox{as } {n \to \infty}.
\end{equation}
We next show  that  $v \in H^1(\rn_+)$ is a weak solution of
\begin{equation}
\left\{
       \begin{array}{rcl}
         -\Delta v + v & = & f(v^2)v \textrm{ in } \rn_+,\\
         v & = & 0 \textrm{ on } \mathbb R ^{N-1}.
       \end{array}
\right. \label{2prn+}
\end{equation}
Effectively, set  $$w_n (x) \doteq \xi\left(\frac{|x|}{\alpha_n} \right)\tilde v_n (x), \quad \forall\, x \in \tilde \Omega_n,
$$
where $\alpha_n > 0$ is given in \eqref{baricentro6} and $\xi \in \cinfsc (\re_+)$ is such that $\xi(t) = 1$, $t \in [0,\frac{1}{2}]$, $\xi(t) = 0$, $t \ge {2}/{3}$. Thus, $w_n \in H_0^1(\tilde\Omega_n)\subset H_0^1(\rn_+)$ and  $w_n (x) \to v(x)$ almost every  $x\in \rn_+$, as $n\to \infty$. Since  $(w_n) \subset H_0^{1}(\rn_+)$ is bounded, there is $w \in H^1_0(\rn_+)$ such that  $w_n \tofraco w$  weakly in   $H^1_0(\rn_+)$. By the Sobolev imbedding theorem,
$w_n (x) \to w(x)$ almost every $x\in \rn_+$, as $n\to \infty$. As the limit is unique, $v \equiv w$ in $H_0^1(\rn_+)$. Taking $\phi \in \cinfsc(\rn_+)$, gives $\supp \phi \subset \tilde \Omega_n$ for every $n$ sufficiently large. By \eqref{baricentro5} and the definition of $\tilde v_n$, we have
\begin{equation}
\label{baricentro9}
\int_{\tilde \Omega_n}( \grad \tilde v_n \grad \phi + \tilde v_n \phi - f(\tilde v_n^2) \tilde v_n \phi )= o_1(n),
\end{equation}
for every $n$ sufficiently large. From \eqref{baricentro7}, after taking $n \to\infty$ in \eqref{baricentro9}, we find
\[
\int_{\rn_+}( \grad v \grad \phi + v \phi - f(v^2)v\phi )= 0.
\]
Since $\phi$ is arbitrary,  the function $v$ is a weak solution of \eqref{2prn+}. Let $J_{\tilde \Omega_n}: H^1(\tilde\Omega_n, \tilde\Gamma_{\lambda_n 0}) \to \re$ be the functional associated with the problem
\begin{equation}
\left\{\begin{array}{rcl}
        -\Delta v + v & = & f(v^2)v \textrm{ in } \tilde\Omega_n,\\
        \displaystyle{\frac{\partial v}{\partial \nu}} & = & 0 \textrm{ on } \tilde\Gamma_{1 \lambda_n}, \\
        v & = & 0 \textrm{ on } \tilde\Gamma_{0 \lambda_n}.\\
       \end{array}
\right. \label{2potl}
\end{equation}
Using that $v$ is a weak solution of \eqref{2prn+}, Fatou lemma and  \eqref{baricentro5}, we have
\begin{eqnarray*}
\ci & = & \liminf_{n \to\infty} J_{\tilde\Omega_n} (\tilde v_n) =\liminf_{n \to\infty} \frac{1}{2} \int_{\tilde \Omega_n}( f(\tilde v_n^2)\tilde v_n^2 - F(\tilde v^2_n) )\\
    & = & \liminf_{n \to\infty} \frac{1}{2} \int_{\rn_+} \chi_{\tilde \Omega_n}( f(\tilde v_n^2)\tilde v_n^2 - F(\tilde v^2_n))
     \ge  \frac{1}{2} \int_{\rn_+} (f(v^2)v^2 - F(v^2) )\\
    & \ge & c_{\rn_+}  \ge  \ci,
\end{eqnarray*}
that is $\ci =  c_{\rn_+}$. However, $\ci =  c_{\rn_+} \ge c_{\rn} = 2 \ci$, which is impossible, and {Claim} I is proved.

\smallskip

\textbf{Claim II.} \emph{Given any\,  $\epsilon > 0$, there exists $R = R(\epsilon) >0$ such that
\begin{equation}\label{claim2}
\displaystyle{\lim\limits_{n \to \infty}\int_{\Omega_{\lambda_n} \cap B_R(y_n)}\left[ f(v_n^2)v_n^2 - F(v_n^2)\right] \ge \ci - \epsilon.}
\end{equation}}

\smallskip
Indeed, we first show that the function $v$ given by \eqref{baricentro7}-\eqref{baricentro8} satisfies  $J_\infty (v) = \ci$ and  $v$ is a solution of \eqref{2pinf}. Consider  $\phi \in C^{\infty}_{0}(\mathbb{R}^{N})$ such that
$\phi = 1$, on  $B_{1}(0)$, $\phi=0$, on $B_{2}^{c}(0)$, $0 \leq \phi \leq 1$, and define
\[
\phi_{T}(x)=\phi\left(\frac{x}{T}\right), \quad \forall x \in \rn, T>0.
\]
Hence, the sequence $\phi_T \tilde{v}_n$ is bounded in $H^{1}(\tilde{\Omega}_n, \tilde{\Gamma}_{0 \lambda_n})$ and  $\phi_T v \to v$ in $H^{1}(\mathbb{R}^{N}_{+})$, as $T\to \infty$. By \eqref{baricentro5}, we have
\[
\int_{\tilde{\Omega}_n}\grad \tilde{v}_n \nabla(\phi_T \tilde{v}_n)+ \int_{\tilde{\Omega}_n} |\tilde{v}_n|^2 \phi_T = \int_{\tilde{\Omega}_n}f(\tilde{v}_n^{2})\tilde{v}_n^{2} \phi_T + o_n(1),
\]
that is,
\begin{equation}
\label{baricentro11}
\int_{\tilde{\Omega}_n}|\nabla \tilde{v}_n|^2 \phi_T + \int_{\tilde{\Omega}_n} \tilde{v}_{n} \nabla \tilde{v}_n \nabla \phi_T + \int_{\tilde{\Omega}_n} |\tilde{v}_n|^2 \phi_T = \int_{\tilde{\Omega}_n}f(\tilde{v}_n^{2})\tilde{v}_n^{2} \phi_T + o_n(1).
\end{equation}
We now proceed to verify that
\begin{eqnarray}
\label{baricentro12}
& &\int_{\tilde{\Omega}_n} \tilde{v}_{n} \nabla \tilde{v}_n \nabla \phi_T \to \int_{\mathbb{R}^{N}_{+}} v  \nabla v \nabla \phi_T, \\
\label{baricentro13}
& &\int_{\tilde{\Omega}_n} |\tilde{v}_n|^2 \phi_T \to \int_{\mathbb{R}^{N}_{+}} |v|^2 \phi_T,\\
& &
\label{baricentro14}
\int_{\tilde{\Omega}_n}f(\tilde{v}_n^{2})\tilde{v}_n^{2} \phi_T \to \int_{\mathbb{R}^{N}_{+}}f(v^2)v^2 \phi_T,
\end{eqnarray}
as $n \to \infty$.  Let $\epsilon > 0$ and $T > 1$ be arbitrary numbers. Fix  $t>0$ to be appropriately chosen and define
\[
E_t \doteq \left\{ x \in B_{2T}(0) ; 0 \le x^N \le t \right\}.
\]
Using that $(\|\tilde v_n\|_{\tilde \Omega_n})_n$ is bounded and
Holder inequality, we obtain
\[
\int_{E_t} |\tilde v_n|^2  \le  \left( \int_{E_t} |\tilde v_n|^{2^*} \right)^{\frac{2}{2^*}} \left( \int_{E_t} |1|^{N} \right)^{\frac{2}{N}}
                           \le  M |E_t|^{\frac{2}{N}} \le M T^{(N-1) \frac{2}{N}} t^{\frac{2}{N}}
\]
and
\[
\int_{E_t} |\tilde v_n|^q  \le  \left( \int_{E_t} |\tilde v_n|^{2^*} \right)^{\frac{q}{2^*}} \left( \int_{E_t} |1|^{\alpha} \right)^{\frac{q}{\alpha}}
                           \le  M |E_t|^{\frac{2}{N}} \le M T^{(N-1) \frac{q}{\alpha}} t^{\frac{q}{\alpha}},
\]
for some positive constant $M$, where $\alpha = 2^*/(2^*-q)$.  Set $\kappa \doteq \max\{{\alpha}/{q}, {N}/{2}\}$ and take $t \doteq \epsilon^{\kappa} T^{1-N}$. Thus,
\[
T^{(N-1) \frac{2}{N}} t^{\frac{2}{N}} = \epsilon^{\kappa\frac{2}{N}} \quad \textrm{and} \quad T^{(N-1) \frac{q}{\alpha}} t^{\frac{q}{\alpha}} = \epsilon^{\kappa\frac{q}{\alpha}},\quad \textrm{with } \min\{\kappa{2}/{N}, \kappa{q}/{\alpha}\} = 1 > 0,
\]
 and consequently
 \begin{equation}
\label{baricentro15}
\lim_{\epsilon \to 0} \epsilon^{\kappa\frac{q}{\alpha}} = \lim_{\epsilon \to 0} \epsilon^{\kappa\frac{2}{N}} = 0.
\end{equation}
By choice of $t$, we have
\begin{equation}
\label{baricentro16}
\int_{E_t} |\tilde v_n|^2 \le \epsilon^{\kappa\frac{2}{N}} M \quad \textrm{and} \quad \int_{E_t} |\tilde v_n|^q \le \epsilon^{\kappa\frac{q}{\alpha}} M
\end{equation}
We observe that by \eqref{baricentro7}-\eqref{baricentro7},  $v$ also satisfies \eqref{baricentro16}. Furthermore, $B_{2T}\backslash E_t \subset \tilde \Omega_n$, provided that $n$ is sufficiently large. Applying Holder inequality, \eqref{baricentro7}, \eqref{baricentro8} and \eqref{baricentro16}, for every $n$ sufficiently large, we get
\begin{eqnarray*}
\lefteqn{\left|\int_{\tilde \Omega_n} \tilde v_n \nabla \phi_T \nabla \tilde v_n - \int_{\mathbb R ^N_+} v \nabla \phi_T \nabla v\right| }\\ &\le  & \left|\int_{B_{2T} \backslash E_t} \tilde v_n \nabla \phi_T \nabla \tilde v_n - v \nabla \phi_T \nabla v\right| + \left|\int_{E_t} \tilde v_n \nabla \phi_T \nabla \tilde v_n\right| + \left|\int_{E_t} v \nabla \phi_T \nabla v\right|  \\
                   & \leq & \left|\int_{B_{2T} \backslash E_t} \nabla \phi_T \nabla \tilde v_n (\tilde v_n - v)\right| + \left|\int_{B_{2T} \backslash E_t} v \nabla \phi_T (\nabla \tilde v_n - \nabla v)\right|  + M\int_{E_t} (|\tilde v_n|^2 + |v|^2) \\
                   & \le &  o_n(1) + {2M  \epsilon^{\kappa\frac{2}{N}}}.
\end{eqnarray*}
From \eqref{baricentro15} and the fact that $\epsilon$ can be chosen arbitrarily small, we obtain that \eqref{baricentro12} holds for every $T>0$. We can proceed analogously to proof of \eqref{baricentro13}. In order to verify \eqref{baricentro14}, we combine  $(f_{6})$ with \eqref{baricentro16}, to obtain
\begin{eqnarray}
\label{baricentro17}
\int_{E_t} f(\tilde{v}_n^{2})\tilde{v}_n^{2}  \le  \epsilon \int_{E_t} \tilde v_n^2 + C_\epsilon \int_{E_t} \tilde v_n^q
                                              \le  M (\epsilon^{\kappa\frac{q}{\alpha}} + \epsilon^{\kappa\frac{2}{N}}).
\end{eqnarray}
From \eqref{baricentro7} and \eqref{baricentro17}, we have
\begin{eqnarray*}
\lefteqn{\left|\int_{\tilde \Omega_n} f(\tilde{v}_n^{2})\tilde{v}_n^{2} \phi_T - \int_{\mathbb R^N_+} f(v^2)v^2 \phi_T \right|}\\
 & \le & \left|\int_{B_{2T}\backslash E_t} (f(\tilde{v}_n^{2})\tilde{v}_n^{2} - f(v^2)v^2) \phi_T \right| + \left|\int_{E_t} f(\tilde{v}_n^{2})\tilde{v}_n^{2}\right| + \left|\int_{E_t} f(v^2)v^2 \right| \\
                              & \le & o_n (1) + M (\epsilon^{\kappa\frac{q}{\alpha}} + \epsilon^{\kappa\frac{2}{N}}).
\end{eqnarray*}
From \eqref{baricentro15} and the fact that $\epsilon$ can be chosen arbitrarily small, we obtain that \eqref{baricentro14} holds for every $T>0$. Combining \eqref{baricentro11}-\eqref{baricentro14} with Fatou lemma, we get
\[
\int_{\mathbb{R}^{N}_{+}}|\nabla v|^2 \phi_T + \int_{\mathbb{R}^{N}_{+}} v \nabla v \nabla \phi_T + \int_{\mathbb{R}^{N}_{+}}|v|^2 \phi_T \leq  \int_{\mathbb{R}^{N}_{+}}f(v^2)v^2 \phi_T,
\]
for every $T>0$. Finally, taking $T \to +\infty$, we find
\begin{equation}
\label{baricentro10}
\int_{\rn_+} {[|\grad v|^2 + v^2]} \le \int_{\rn_+} f(v^2)v^2.
\end{equation}
From $(f_1)-(f_{4})$, there exists $t_0 > 0$ such that $t_0 v \in N_{\infty}$. By \eqref{baricentro10}, we have $0 < t_0 \le 1$. Suppose that  $t_0 < 1$. In this case, using that the function $s \to f(s)s-F(s)$ is increasing in $[0,+\infty)$, by $(f_4)$,  Fatou lemma and \eqref{baricentro8}, we have
\begin{eqnarray*}
\ci & = & \liminf_{n \to \infty}[ J_{\tilde\Omega_n} (\tilde v_n) - \frac{1}{2} J'_{\tilde \Omega_n} (\tilde v_n)\tilde v_n]  = \liminf_{n \to \infty}\frac{1}{2}\int_{\tilde \Omega_n}[f(\tilde v_n^2)\tilde v_n^2 - F(\tilde v_n^2)] \\
    & \ge & \frac{1}{2}\int_{\rn_+}[f(v^2)v^2 - F(v^2)]
     >  \frac{1}{2}\int_{\rn_+} [f(t^2_0 v^2) t_0^2 v^2 - F(t_0^2 v^2) ] \\
    & = & { J_\infty(t_0 v) - \frac{1}{2} J'_{\infty} (t_0 v) t_0 v \ge \ci,}
\end{eqnarray*}
which is impossible. Hence, $t_0 = 1$, and consequently  $v \in N_{\infty}$. Furthermore, $v$ satisfies
\begin{equation}
\label{baricentro19}
\ci \ge \frac{1}{2}\int_{\rn_+}[ f(v^2)v^2 - F(v^2)] = J_\infty(v) -\frac{1}{2}J'_\infty(v)v = J_\infty(v) \ge \ci.
\end{equation}
We conclude that $J_\infty (v) = \ci$ and  $v$ is a solution of \eqref{2pinf}.  By  (\ref{baricentro19}), given any $\epsilon > 0$, there exists $R > 0$ such that
\[
\frac{1}{2}\int_{\rn_+ \cap B_R(0)} [f(v^2)v^2 - F(v^2)] \ge \ci - \epsilon.
\]
Since $\chi_{B_R \cap \tilde \Omega_n} (x) \tilde v_n(x) \to \chi_{B_R^+}(x) v(x)$ almost every  $x\in \rn_+$, as $n\to \infty$, by
Fatou lemma we have
\begin{eqnarray*}
\liminf_{n \to \infty} \frac{1}{2}\int_{B_R (y_n) \cap \Omega_n}[ f(v_n^2) v_n^2 - F(v_n^2)] & = & \liminf_{n \to \infty} \frac{1}{2}\int_{B_R \cap \tilde \Omega_n}[ f(\tilde v_n^2)\tilde v_n^2 - F(\tilde v_n^2)] \\
         & \ge & \frac{1}{2}\int_{\rn_+ \cap B_R(0)} [f(v^2)v^2 - F(v^2)] \\
         & \ge & \ci - \epsilon,
\end{eqnarray*}
which completes the proof of Claim II.

 \medskip
 We are now ready to show \eqref{baricentro2}. By \eqref{baricentro4} and the Sobolev embedding theorem, the sequences {$\{t_n |u_n | \} \subset  H^1(\Omega_n, \Gamma_{0 \lambda_n})$ and   $\{v_n \} \subset H^1(\Omega_n, \Gamma_{0 \lambda_n}) $} have the same limit. Hence, Claim II is also valid for  $\{t_n |u_n|\}_n$, that is,
\[
\liminf_{n \to \infty} \frac{1}{2}\int_{B_R(y_n) \cap \Omega_{\lambda_n}} [f(|t_n u_n |^2) |t_n u_n|^2 - F(|t_n u_n|^2)] \ge \ci - \epsilon.
\]
From this, \eqref{baricentro33} and  $(f_{5})$, we have
\begin{equation}
\label{baricentro21}
\liminf_{n \to \infty}\int_{\Omega_{\lambda_n}\backslash B_R (y_n)}  C |t_n u_n|^q \le \epsilon.
\end{equation}
By Claim I, we can assume that  $y_n \in \Gamma_{1 \lambda_n}$, i.e. ${y_n}/{\lambda_n} \in \Gamma_1$ and ${y_n}/{\lambda_n} \to x_0 \in \overline\Gamma_1$, as $ n\to \infty$, because $\overline \Gamma_1$ is a compact set.   Take $j \in \{1, \dots, N\}$. From the definition of the barycenter, we have
\begin{eqnarray*}
\left| \frac{\beta_{\lambda_n}^j(u_n)}{\lambda_n} - x_0^j \right| & \le & \frac{\displaystyle\int_{\Omega_{\lambda_n}} \left|\frac{x^j}{\lambda_n} - x_0^j\right||t_n u_n|^q}{\displaystyle\int_{\Omega_{\lambda_n}}|t_n u_n|^q}.
\end{eqnarray*}
Using Lemma \ref{lionsnosso} and the fact that $t_n |u_n| \in M_{\lambda_n}$,  we may assume that
\[
\int_{\Omega_{\lambda_n}}|t_n u_n|^q \ge \gamma > 0, \quad \forall n \in \mathbb N.
\]
As a consequence,
\begin{eqnarray*}
\gamma \left|\frac{\beta_{\lambda_n}^j(u_n)}{\lambda_n} - x_0^j\right| & \le & \int_{\Omega_{\lambda_n}} \left|\frac{x^j}{\lambda_n} - x_0^j\right||t_n u_n|^q \\
                       & = & \int_{\Omega_{\lambda_n} \cap B_R (y_n)} \left|\frac{x^j}{\lambda_n} - x_0^j\right| |t_n u_n|^q + \int_{\Omega_{\lambda_n} \backslash B_R (y_n)} \left|\frac{x^j}{\lambda_n} - x_0^j\right||t_n u_n|^q \\
                       & \le & \int_{\Omega_{\lambda_n} \cap B_R (y_n)} \left|\frac{x^j}{\lambda_n} - \frac{y_n^j}{\lambda_n}\right||t_n u_n|^q + \int_{\Omega_{\lambda_n} \cap B_R (y_n)} \left|\frac{y_n^j}{\lambda_n} - x_0^j\right||t_n u_n|^q \\
                       &     & + \int_{\Omega_{\lambda_n} \backslash B_R(y_n)} \left|\frac{x^j}{\lambda_n} - x_0^j\right||t_n u_n|^q \\
                       & \le & \frac{R}{\lambda_n} \int_{\Omega_{\lambda_n}} |t_n u_n|^q + \left|\frac{y_n}{\lambda_n} - x_0\right| \int_{\Omega_{\lambda_n}} |t_n u_n|^q  \\
                       &  & + \textrm{ diam}(\Omega) \int_{\Omega_{\lambda_n} \backslash B_R(y_n)} |t_n u_n|^q \\
                       & = & \left(\frac{R}{\lambda_n} + \left|\frac{y_n}{\lambda_n} - x_0\right|\right) \int_{\Omega_n} |t_n u_n|^q + \textrm{diam}(\Omega) \int_{\Omega_{\lambda_n} \backslash B_R (y_n)} |t_n u_n|^q.
\end{eqnarray*}
From \eqref{baricentro21} and the fact that the sequence $(\| t_n u_n \|_{A_{\lambda_n}})_n$ is bounded and ${y_n}/{\lambda_n} \to x_0$, we find
\[
0 \le \liminf_{n \to \infty}  \left|\frac{\beta_{\lambda_n}^j(u_n)}{\lambda_n} - x_0^j\right| \le  \textrm{diam}(\Omega) \frac{\epsilon}{\gamma C}, \quad \forall\, j \in \{1, \dots, N\}.
\]
Since $\epsilon >0$ is arbitrary, we can find a subsequence (not renamed) such that
\[
\textrm{dist}\left(\frac{\beta_{\lambda_n} (u_n)}{\lambda_n}, \Gamma_1\right) \to  0,\quad \mbox{as $n\to \infty$}.
\]
We conclude that $\textrm{dist}(\beta_{\lambda_n}(u_n), \Gamma_{1 \lambda_n}) \le \lambda_n r$, for every $n$ sufficiently large, hence that \eqref{baricentro2} holds, and the proposition follows. \cqd

{Taking $\epsilon^*>0$  given by Proposition \ref{baricentro},  we define   $b^*_\lambda = b_{\lambda} + \epsilon^*$. As a consequence of Propositions  \ref{phi}, \ref{bltoci} and \ref{baricentro}, we obtain the following result which  is the key point in the comparison of the
topology of  the sublevel sets of the functional $I_\lambda$ with
 that of $\Gamma_{1 \lambda}$.}
\begin{lemma} \label{phibar} There exists  $\lambda^*>0$  such that
\[
\Phi_{\lambda}((\Gamma^-_1)_{\lambda}) \subset M^{b^*_\lambda}_{\lambda} \quad \mbox{and}\quad \beta_\lambda(M_\lambda^{b^*_\lambda}) \subset (\Gamma_1^+)_{\lambda },
\]
for every  $\lambda > \lambda^*$, where $M_\lambda^{b^*_\lambda} \doteq I_{\lambda}^{b^*_\lambda} \cap M_{\lambda}$.
\end{lemma}

\dem\ By Proposition \ref{baricentro}, there exists $\lambda_1> 0$ such that $\beta_\lambda(M_\lambda^{b^*_\lambda}) \subset (\Gamma_1^+)_{\lambda }$.  From Propositions  \ref{phi} and \ref{bltoci}, we have
\begin{equation}
\label{phibar1}
\lim_{\lambda \to \infty} (I_{\lambda}(\Phi_\lambda (y)) - b_{\lambda}) = 0,
\end{equation}
independent of $y\in \Gamma^-_{1\lambda}$. Thus, for this $\epsilon^* > 0$ there exits  $\lambda_2= \lambda_2(\epsilon^*)  > 0$  such that
\[
I_{\lambda} (\Phi_\lambda(y)) \le b_{\lambda} + \epsilon^*,
\]
for every $\lambda > \lambda_2$ and $y \in (\Gamma^-_1)_{\lambda}$. Set
$\lambda^* \doteq \max\{\lambda_1, \lambda_2\}$. Hence,
\[
\Phi_\lambda((\Gamma^-_1)_{\lambda}) \subset M_\lambda^{b^*_\lambda} \quad \textrm{and} \quad \beta_\lambda(M_\lambda^{b^*_\lambda}) \subset (\Gamma^+_1)_{\lambda},\quad \forall\, \lambda > \lambda^*.
\]
\cqd

\section{Proof of Theorem \ref{2t1}}
We begin by stating a comparison of the topology of  the sublevel $M_\lambda^{b^*_\lambda}$ with that of $\Gamma_{1 \lambda}$.
\begin{lemma} \label{catml} Let  $\lambda^* > 0$ be as in Lemma \ref{phibar}. Then,
\[
cat_{M_{\lambda}^{b^*_{\lambda}}}(M_{\lambda}^{b^*_{\lambda}}) \ge cat_{\Gamma_{1\lambda}}(\Gamma_{1\lambda}),
\]
for every  $\lambda > \lambda^*$.
\end{lemma}

\dem \  The proof proceeds along the same lines as the proof of  \cite[Lemma 3.5]{benci-cerami}.   Suppose that  $cat_{M_{\lambda}^{b^*_{\lambda}}}(M_{\lambda}^{b^*_{\lambda}}) = m$. Thus, $M_{\lambda}^{b^*_{\lambda}} = \Upsilon_1 \cup \dots \cup \Upsilon_m$, where $\Upsilon_j$ is closed and contractible in $M_{\lambda}^{b^*_{\lambda}}$, for  $j = 1, \dots, m$. Hence, there exists  $h_j \in C([0,1] \times \Upsilon_j, M_{\lambda}^{b^*_{\lambda}})$ such that
$h_j(0,u) = u$,  $h_j(1,u) = u_j \in M_{\lambda}^{b^*_{\lambda}}$ for every $u \in \Upsilon_j$ and $j = 1, \dots, m$, for some   $u_j \in M_{\lambda}^{b^*_{\lambda}}$ fixed. Set  $B_j := \Phi_\lambda^{-1}(\Upsilon_j)$, $j = 1, \dots, m$, which are closed in
$\Gamma_{1 \lambda}^-$. By  Proposition \ref{phibar}, we have
\[
\Gamma_{1 \lambda}^- = \bigcup^m_{j = 1} B_j.
\]
Using Proposition \ref{phibar} again, the maps $g_j : [0,1] \times B_j \to \Gamma_{1 \lambda}^+$ given by
\[
g_j (t,y) := \beta_\lambda(h_j(t,\Phi_\lambda(y))), \quad \forall j \in \{1, \ldots, m\},
\]
are well defined. In addition, $g_j \in C([0,1] \times B_j, \Gamma_{1 \lambda}^+)$ and
\[
g_j(0,y) = y, \, g_j(1,y) = y_j \in \Gamma_{1 \lambda}^+, \, \textrm{ for every } y \in B_j, \, j = 1, \dots, m,
\]
{and $y_j \in \Gamma_{1 \lambda}^+$ fixed, and so $cat_{\Gamma_{1 \lambda}^+} \Gamma_{1 \lambda}^- \le m$. Recalling that $\Gamma_{1 \lambda}^+$ and $\Gamma_{1 \lambda}^-$ are homotopically equivalent to $\Gamma_{1 \lambda}$, it follows that $cat_{\Gamma_{1 \lambda}} \Gamma_{1 \lambda} = cat_{\Gamma_{1 \lambda}^+} \Gamma_{1 \lambda}^-$, and hence $cat_{\Gamma_{1 \lambda}} \Gamma_{1 \lambda} \leq m$, which completes the proof.  \cqd }

\textbf{Proof of Theorem \ref{2t1}}.
{Take $\epsilon^*>0$ given by  Proposition \ref{baricentro}, $\lambda^* > 0$ given by  Proposition \ref{phibar}, and  suppose $\lambda \ge \lambda^*$.
If $b^*_\lambda =     b_\lambda + \epsilon$ is a critical value for every $\epsilon \in (0, \epsilon^*]$ then $I_\lambda$ has infinitely many critical values and the proof is complete.  Otherwise, we can assume  that $b^*_\lambda$  is a regular value of $I_\lambda$.
  Since  $M_{\lambda}^{b_\lambda^*}$ is a closed set in $M_\lambda$, by Proposition \ref{psr}, the restriction of $I_\lambda$ to  $M_{\lambda}^{b_\lambda^*}$ satisfies the  $(PS)_d$ condition for every   $d \in \re$. Hence, by  the Ljusternik-Schnirelman theory
and  Lemma  \ref{catml}, we obtain $cat_{\Gamma_{1 \lambda}}(\Gamma_{1 \lambda})$ critical points of  $I_\lambda \big|_{M_{\lambda}^{b_\lambda^*}}$. By Corollary \ref{pcr}, each of these critical points is a critical point of $I_\lambda$. }
\cqd

\section{Morse theory for  $I_\lambda$}
In this section we see how the homology groups of the sets $\Gamma_{1\lambda}, (\Gamma_{1}^-)_\lambda, (\Gamma_{1}^+)_\lambda$ and $M_{\lambda}^{b^*_\lambda}$
are related.  For the convenience of the reader, we repeat  the relevant material from \cite[Section 5]{benci-cerami} adapted to our case, thus making the exposition self-contained.
\begin{lemma}\label{pmlpg} Let  $\lambda^* > 0$ be as in Lemma \ref{phibar}. Then,
\begin{equation*}
\mathcal P_t (M_{\lambda}^{b^*_{\lambda}}) = \mathcal P_t (\Gamma_{1 \lambda}) + \mathcal Q(t),
\end{equation*}
for every  $\lambda \ge \lambda^*$, where $\mathcal Q$ is a polynomial with non-negative coefficients.
\end{lemma}

\dem\  Setting $\lambda \ge \lambda^*$, the function $\Phi_\lambda : (\Gamma_{1}^-)_\lambda \to  M_{\lambda}$ given by \eqref{phi} induces the homomorphism
$(\Phi_\lambda)_k : H_k(\Gamma^-_{1 \lambda}) \to H_k(M_\lambda^{b_\lambda^*})$
between the $k$-th homology groups. Since $\Phi_\lambda$  is a {{injective}} function, so also is $(\Phi_\lambda)_k$. Hence, $\textrm{dim}H_k(\Gamma^-_{1 \lambda}) \geq \textrm{dim}H_k(M_\lambda^{b_\lambda^*})$, and the result follows from the definition of the Poincar\'{e} polynomials and the fact that $\Gamma_{1 \lambda}^-$ and $\Gamma_{1 \lambda}$ are homotopically equivalent.
\cqd
\begin{lemma} \label{pibidpml} Let  $\lambda^* > 0$ be as in Lemma \ref{phibar}, $\lambda \ge \lambda^*$, $\delta \in (0, \delta_0)$, for $\delta_0$ given by  Proposition \ref{limnehari}, and  $b \in \left.(\delta, \infty\right.]$ a noncritical level of $I_\lambda$. Then,
\begin{equation*}
\mathcal P_t (I^b_\lambda, I^\delta_\lambda) = t \mathcal P_t(M^b_\lambda)
\end{equation*}
\end{lemma}

\dem \   The proof proceeds along the same lines as the proof of \cite[Lemma 5.2]{benci-cerami}. \cqd

\begin{lemma} \label{phalibpg} Let $\lambda^*$, $\lambda$ and $\delta$ be as in Lemma \ref{pibidpml}. Then
\begin{equation}
\label{ph1}
\mathcal P_t (I^{b^*_\lambda}_\lambda, I_\lambda^\delta) = t \mathcal P_t (\Gamma_{1\lambda}) + t\mathcal Q (t)
\end{equation}
and
\begin{equation}
\label{hdoispol2}
\mathcal P_t (H^1_\al (\Omega_\lambda, \Gamma_{0 \lambda}), I_\lambda^\delta) = t \mathcal P_t (M_\lambda)=t,
\end{equation}
where $\mathcal Q$ is a polynomial with non-negative coefficients.
\end{lemma}

\dem\   As in the proof of Theorem  \ref{2t1}, we can assume that  $b_\lambda^*$ is a regular value. Applying Lemma \ref{pibidpml}, for  $b = b_\lambda^*$, and Lemma \ref{pmlpg},  we get  \eqref{ph1}. Using that $M_{\lambda}$ is homeomorphic to the unit sphere in $H^1_{A_\lambda} (\Omega_\lambda, \Gamma_{0\lambda})$,
which is contractible (see \cite[Example 1B.3]{ah}), we have that $M_{\lambda}$ is
 contractible. Hence, ${\rm{dim}}H^k(M_\lambda) = 1$ if $k = 0$ and ${\rm{dim}} H^k(M_\lambda) = 0$ if $k \neq  0$. Finally, \eqref{hdoispol2} is obtained by again invoking Lemma \ref{pibidpml}, for  $b = \infty$.

\begin{lemma} \label{hphalibpg} Let $\lambda^*$, $\lambda$ and $\delta$ be as in Lemma \ref{pibidpml}. Then
\begin{equation}\label{ph2}
\mathcal P_t (H^1_\al(\Omega_\lambda, \Gamma_{0 \lambda}), I_\lambda^{b^*_\lambda}) = t^2 [ P_t (\Gamma_{1 \lambda}) + \mathcal Q(t) - 1 ],
\end{equation}
where $\mathcal Q$ is a polynomial with non-negative coefficients.
\end{lemma}

\dem\  We follow Benci and Cerami \cite{benci-cerami}  in considering the exact sequence:
\begin{align*}
\dots \longrightarrow H_k (H_\al^1(\Omega_\lambda, \Gamma_{0 \lambda}), I^\delta_\lambda) \stackrel{j_k}\longrightarrow & H_k (H_\al^1(\Omega_\lambda, \Gamma_{0 \lambda}), I^{b^*_\lambda}_\lambda) \stackrel{\partial_k}\longrightarrow \\
\stackrel{\partial_k}\longrightarrow H_{k-1} (I^{b^*_\lambda}_\lambda, I^\delta_\lambda) & \stackrel{i_{k-1}}\longrightarrow H_{k-1} (H_\al^1(\Omega_\lambda, \Gamma_{0 \lambda}), I^\delta_\lambda) \longrightarrow \dots
\end{align*}
From \eqref{hdoispol2}, we obtain
$
\textrm{dim} H_k(H_\al^1(\Omega_\lambda, \Gamma_{0 \lambda}), I^\delta_\lambda) = 0, \, \forall k \neq 1.
$
If we combine this with the fact that the sequence is exact , we see that  $\partial_k$ is a  isomorphism for every $k \ge 3$. Hence,
\begin{equation}
\label{hphalibpg1}
\textrm{dim} H_k(H_\al^1(\Omega_\lambda, \Gamma_{0 \lambda}), I^{b^*_\lambda}_\lambda) = \textrm{dim} H_{k-1} (I^{b^*_\lambda}_\lambda, I^\delta_\lambda), \, \forall k \ge 3.
\end{equation}
For $k = 2$, we have
\begin{align*}
\dots \longrightarrow H_2 (H_\al^1(\Omega_\lambda, \Gamma_{0 \lambda}), I^\delta_\lambda) \stackrel{j_2}\longrightarrow & H_2 (H_\al^1(\Omega_\lambda, \Gamma_{0 \lambda}), I^{b^*_\lambda}_\lambda) \stackrel{\partial_2}\longrightarrow  \\
\stackrel{\partial_2}\longrightarrow H_{1} (I^{b^*_\lambda}_\lambda, I^\delta_\lambda) & \stackrel{i_1}\longrightarrow H_{1} (H_\al^1(\Omega_\lambda, \Gamma_{0 \lambda}), I^\delta_\lambda) \longrightarrow \dots
\end{align*}
Since $j_2$ is sobrejective ($j_2$ is the homomorphism induced by the canonic projection) and $\textrm{dim} H_2 (H_\al^1(\Omega_\lambda, \Gamma_{0 \lambda}), I^\delta_\lambda) = 0$,  by  \eqref{hdoispol2}, we have
\begin{equation}
\label{phalibpg2}
H_2 (H_\al^1(\Omega_\lambda, \Gamma_{0 \lambda}), I_\lambda^{b^*_\lambda}) = j_2 ( H_2 (H_\al^1(\Omega_\lambda, \Gamma_{0 \lambda}), I^\delta_{\lambda}) ) = \{0\}.
\end{equation}
For  $k = 1$,
\begin{align*}
\dots \longrightarrow H_1 (I^{b^*_\lambda}_\lambda, I^\delta_\lambda) \stackrel{i_1}\longrightarrow & H_1 (H_\al^1(\Omega_\lambda, \Gamma_{0 \lambda}), I^\delta_\lambda) \stackrel{j_1}\longrightarrow \\ \stackrel{j_1}\longrightarrow H_1 (H_\al^1(\Omega_\lambda, \Gamma_{0 \lambda}), I^{b^*_\lambda}_\lambda) & \stackrel{\partial_1}\longrightarrow H_0 (I^{b^*_\lambda}_\lambda, I^\delta_\lambda) \longrightarrow \dots
\end{align*}
Using that  $H_\al^1(\Omega_\lambda, \Gamma_{0 \lambda})$ is a connected set, we have
\begin{equation}
\label{hphalibpg3}
H_0 (H_\al^1(\Omega_\lambda, \Gamma_{0 \lambda}), I^{b^*_\lambda}_\lambda) = 0.
\end{equation}
We now claim that $i_1$ is a isomorphism. Indeed, as  $\Gamma_{1 \lambda}     \neq \emptyset$ and  ${\textrm{dim}} H_0 (\Gamma_{1\lambda})$ is the number of connected components of  the set  $\Gamma_{1\lambda}$, we have
$H_0 (\Gamma_{1 \lambda}) \neq \{0\}$. By \eqref{ph1},   $H_1 (I^{b^*_\lambda}_\lambda, I^\delta_\lambda) \neq \{0\}$. From  \eqref{hdoispol2}, we obtain $\textrm{dim} H_1 (H_\al^1(\Omega_\lambda, \Gamma_{0 \lambda}), I^\delta_\lambda) = 1$. Using that
$i_1$ \'{e} is injective, we have $\textrm{dim} H_1 (I^{b^*_\lambda}_\lambda, I^\delta_\lambda) = 1$, and so $i_1$ is a isomorphism. Using that  $i_1$ is a isomorphism and  $j_1$ is sobrejective, we get
\begin{equation}
\label{hphalibpg6}
\textrm{dim} H_1 (H^1_{A_\al}(\Omega_\lambda, \Gamma_{0\lambda}), I^{b^*_\lambda}_\lambda) = 0.
\end{equation}
Combining Lemma \ref{phalibpg} with  \eqref{hphalibpg1} - \eqref{hphalibpg6}, we have
\begin{align*}
\mathcal P_t (H_\al^1(\Omega_\lambda, \Gamma_{0 \lambda}), I^{b^*_\lambda}_\lambda) = & \sum_{k \ge 3} t^k \textrm{dim} H_k (H_\al^1(\Omega_\lambda, \Gamma_{0 \lambda}), I^{b^*_\lambda}_\lambda)\\
=& \sum_{k \ge 3} t^k \textrm{dim} H_{k-1} (I^{b^*_\lambda}_\lambda, I^\delta_\lambda) = t\sum_{k \ge 3} t^{k-1} \textrm{dim} H_{k-1} (I^{b^*_\lambda}_\lambda, I^\delta_\lambda) \\
=& t \left[ P_t (I^{b^*_\lambda}_\lambda, I^\delta_\lambda) - t \, \textrm{dim}H_1 (I^{b^*_\lambda}_\lambda, I^\delta_\lambda) - \textrm{dim}H_0 (I^{b^*_\lambda}_\lambda, I^\delta_\lambda) \right]\\
                                     = & \, t^2 \left[ \mathcal P_t (\Gamma_{1\lambda}) + \mathcal Q(t) - 1 \right].
\end{align*} \cqd
\begin{lemma}\label{c1c2}   Let $\lambda^*$, $\lambda$ and $\delta$ be as in Lemma \ref{pibidpml}.   Suppose that the set  $\mathcal K$ of nontrivial solutions of problem \eqref{2pal} is discrete. Then,
\begin{equation}
\label{c1c21}
\sum_{u \in \, \mathcal C_1} i_t (u) = t \mathcal P_t(\Gamma_{1 \lambda}) + t \mathcal Q(t) + (1+t)\mathcal Q_1(t)
\end{equation}
and
\begin{equation}
\label{c1c22}
\sum_{u \in \, \mathcal C_2} i_t (u) = t^2 [\mathcal P_t(\Gamma_{1\lambda}) + \mathcal Q(t) -1 ] + (1+t)\mathcal Q_2(t),
\end{equation}
where
\[
\mathcal C_1 \doteq \left\{ u \in \mathcal K;\  \delta < I_{\lambda}(u) \le b^*_\lambda \right\}\quad
\mbox{and}\quad
\mathcal C_2 \doteq \left\{ u \in \mathcal K;\  b^*_\lambda < I_{\lambda}(u)\right\},
\]
and $\mathcal Q_i$, $i=1, 2$, is a polynomial with non-negative coefficients.
\end{lemma}
\dem\   Using that $I_\lambda$ satisfies $(PS)$ condition and applying \cite[Theorem 4.3]{chang}, there exists a polynomial $\mathcal Q_1$ with non-negative coefficients such that
\begin{equation}
\sum_{u \in \mathcal C _1} i_t (u) = \mathcal P _t (I^{b^*_\lambda}_\lambda, I^\delta_\lambda) + (1+t)\mathcal Q_1(t).
\label{ph3}
\end{equation}
Hence, \eqref{c1c21} is a consequence of \eqref{ph1} and \eqref{c1c22} follows from \eqref{ph2}.\cqd

\textbf{Proof of Theorem \ref{2t2}.}\ Let $\lambda^*$, $\lambda$ and $\delta$ be as in Lemma \ref{pibidpml}.   Since $I_\lambda$ does not have nontrivial solution below the level $\delta_0$, we have $\mathcal K = \mathcal C_1 + \mathcal C_2$, for $\mathcal C_1$ and $\mathcal C_2$ as in Lemma \ref{c1c2}. Hence,
\[
\sum_{u \in \mathcal K} i_t(u) = \sum_{u \in \mathcal C_1} i_t(u) + \sum_{u \in \mathcal C_2} i_t(u),
\]
Using Lemma \ref{c1c2}, we conclude the proof. \cqd

\textbf{Proof of Corollary \ref{2c1}.}\ This is a direct consequence of Theorem \ref{2t2} and the fact that
$i_t(u) = t^{\mu(u)}$ in the non-degenerate case.  \cqd


\end{document}